\documentclass[11pt,a4paper]{amsart}
\usepackage{amssymb,amsmath,amsthm,mathtools,mathrsfs}
\usepackage{amscd}
\usepackage[margin=25mm]{geometry}
\usepackage[shortlabels]{enumitem}
\usepackage{microtype}
\usepackage{needspace}
\usepackage[dvipsnames]{xcolor}
\usepackage[normalem]{ulem}
\usepackage[hidelinks]{hyperref}
\usepackage{cite}
\usepackage[nodisplayskipstretch]{setspace}
\hypersetup{pdftitle={Second homotopy modules and infinite homotopy stable classes of four-manifolds},pdfauthor={Tommy Hofmann and John Nicholson}}
\newcommand{\Z}{\mathbb Z}
\newcommand{\F}{\mathbb F}
\newcommand{\Q}{\mathbb Q}
\newcommand{\HH}{\widehat H}
\newcommand{\hmod}{\mathcal H}
\DeclareMathOperator{\Aut}{Aut}
\DeclareMathOperator{\End}{End}
\DeclareMathOperator{\nrd}{nrd}
\DeclareMathOperator{\Hom}{Hom}
\DeclareMathOperator{\Ext}{Ext}
\DeclareMathOperator{\Tor}{Tor}
\DeclareMathOperator{\im}{im}
\DeclareMathOperator{\GL}{GL}
\DeclareMathOperator{\SL}{SL}
\DeclareMathOperator{\res}{res}
\DeclareMathOperator{\cor}{cor}
\DeclareMathOperator{\Met}{Met}
\newcommand{\wt}{\widetilde}
\newcommand{\eps}{\varepsilon}
\newcommand{\id}{\mathrm{id}}
\usepackage[capitalize,noabbrev]{cleveref}
\newtheorem{thm}{Theorem}[section]
\newtheorem{prop}[thm]{Proposition}
\newtheorem{lemma}[thm]{Lemma}

\newtheorem{corollary}[thm]{Corollary}
\newtheorem{thmx}{Theorem}
\newtheorem{corx}[thmx]{Corollary}

\theoremstyle{definition}

\theoremstyle{remark}
\newtheorem{remark}[thm]{Remark}

\newenvironment{clist}[1]{\begin{enumerate}[\normalfont #1]}{\end{enumerate}}

\makeatletter
\@namedef{subjclassname@2020}{\textup{2020} Mathematics Subject Classification}
\makeatother

\crefname{thm}{Theorem}{Theorems}
\Crefname{thm}{Theorem}{Theorems}
\crefname{thmx}{Theorem}{Theorems}
\Crefname{thmx}{Theorem}{Theorems}
\crefname{corx}{Corollary}{Corollaries}
\Crefname{corx}{Corollary}{Corollaries}
\crefname{prop}{Proposition}{Propositions}
\Crefname{prop}{Proposition}{Propositions}

\newcommand{\sharedtheoremref}[4]{%
  \AddToHook{env/#1/begin}{\crefalias{#2}{#1}}%
  \crefname{#1}{#3}{#4}%
  \Crefname{#1}{#3}{#4}}
\sharedtheoremref{lemma}{thm}{Lemma}{Lemmas}
\sharedtheoremref{prop}{thm}{Proposition}{Propositions}
\sharedtheoremref{remark}{thm}{Remark}{Remarks}
\sharedtheoremref{corollary}{thm}{Corollary}{Corollaries}
\sharedtheoremref{corx}{thmx}{Corollary}{Corollaries}

\begin{document}
\title[Second homotopy modules and infinite stable classes]{Second homotopy modules and infinite homotopy \\ stable classes of four-manifolds}

\author{Tommy Hofmann}
\address{Naturwissenschaftlich-Technische Fakult\"at, Universit\"at Siegen, Walter-Flex-Strasse 3, 57068 Siegen, Germany}
\email{tommy.hofmann@uni-siegen.de}
\author{John Nicholson}
\address{Max Planck Institute for Mathematics, Vivatsgasse 7, 53111 Bonn, Germany}
\email{math@johnnynicholson.co.uk}
\date{}

\begin{abstract}
We construct infinitely many closed oriented smooth $4$-manifolds which are all stably diffeomorphic but pairwise not homotopy equivalent.
Our examples are distinguished by their second homotopy modules, and provide the first examples which demonstrate that stable diffeomorphism does not determine the second homotopy module. 
They have fundamental group $(\Z/4)^4*\Z$ and arise as doubles of presentation $2$-complexes.
In contrast, we give criteria for isomorphism of second homotopy modules of doubles of $2$-complexes, covering finite abelian, quaternion, dihedral and symmetric fundamental groups.
We also establish non-injectivity of the doubling construction by showing that there exist homotopically distinct presentation $2$-complexes over $Q_{28}$ whose doubles are homotopy equivalent.
\end{abstract}

\hspace{-7mm}

\maketitle

\section{Introduction}

Two closed smooth $4$-manifolds $M$ and $N$ are \textit{stably diffeomorphic} if
\[
 M\#k(S^2\times S^2)\cong N\#k(S^2\times S^2)
\]
for some $k\ge0$.
Kreck--Schafer constructed pairs of closed smooth $4$-manifolds which are stably diffeomorphic but not homotopy equivalent \cite{KS84}.
Kasprowski--Powell--Ray asked whether there exist arbitrarily large families with these properties, or even an infinite family \cite[Question 1.5]{KPR22}.
Recently, Hambleton--Nicholson constructed arbitrarily large finite families \cite[Theorem C]{HN26}, leaving the infinite case open.
They also asked whether stably diffeomorphic closed smooth $4$-manifolds can be distinguished by their second homotopy modules \cite[Question 1.5]{HN26}.

Our main result is the following, which answers both questions in the affirmative.

\begin{thmx}\label{thm:main}
Let $G=(\Z/4)^4*\Z$. There exist infinitely many closed oriented smooth $4$-manifolds $M_1,M_2,\ldots$ with fundamental group $G$ such that
\[
M_i\#(S^2\times S^2)\cong M_j\#(S^2\times S^2) \quad \text{for all }\, i,j
\]
by orientation-preserving diffeomorphisms. Their second homotopy modules $\pi_2(M_i)$ are pairwise non-isomorphic as $\Z G$-modules, even modulo the action of $\Aut(G)$. In particular, the manifolds $M_i$ are pairwise not homotopy equivalent.
\end{thmx}

Two $\Z G$-modules are \textit{$\Aut(G)$-isomorphic} if they become isomorphic after twisting one action by an automorphism of $G$.
We do not require homotopy equivalences to preserve orientation. The \textit{homotopy stable class} of $M$ is the set $\mathcal{S}^{\text{st}}_h(M)$ of homotopy types of manifolds stably diffeomorphic to $M$, as in Conway--Crowley--Powell--Sixt \cite{CCPS24}. Thus \cref{thm:main} gives an infinite homotopy stable class, distinguished already by second homotopy modules.

Infinite homotopy stable classes in dimensions $4k\ge8$ were constructed by Conway--Crowley--Powell--Sixt \cite[Theorem 1.1]{CCPS24}. In dimension $4$, Conway--Crowley--Powell obtained an infinite family of stably homeomorphic topological manifolds with nonempty boundary, distinguished up to homotopy equivalence of pairs \cite[Theorem 1.1]{CCP23}. 

As in the work of Kreck--Schafer \cite{KS84} and Hambleton-Nicholson \cite{HN26}, 
our manifolds are constructed by associating to a finite $2$-complex $X$ the \textit{double} $M(X) = V \cup_{\partial V} (-V)$ of a smooth $4$-dimensional handle thickening $V \simeq X$. 
For $n \ge 0$, let $X_n$ be the presentation $2$-complex of
\[ 
\mathcal P_n=\langle x,y,z,w,t\mid\mathcal R, \ [t^n,x^8]x^4\rangle
 \quad \text{for $n \ge 0$}
\]
where $\mathcal R=\{ x^4y^{-4},\ [x,y]x^4,\ z^4,\ w^4,\ [x,z],\ [x,w],\ [y,z],\ [y,w],\ [z,w]\}$. Each presentation presents $G=(\Z/4)^4*\Z$ such that $(\Z/4)^4 = \langle x,y,z,w\rangle$ and $\Z = \langle t \rangle$ (see \cref{lem:groups}).
The manifolds in \cref{thm:main} are $M_i=M(X_{n_i})$, where the $n_i$ are chosen so that the second homotopy modules lie in distinct $\Aut(G)$-orbits (see Proposition~\ref{prop:unmarked}).

\subsection{Cancellation for finite fundamental groups}

Hambleton--Kreck showed that there cannot exist infinitely many stably diffeomorphic closed smooth $4$-manifolds with finite fundamental group which are homotopically distinct \cite[p87]{HK88}.
If $M$ and $N$ are stably diffeomorphic with fundamental group $G$, then 
\[ \pi_2(M) \oplus \Z G^{2k} \cong \pi_2(N) \oplus \Z G^{2k} \] 
for some $k \ge 0$, after identifying the fundamental groups by the stable diffeomorphism. It follows that $\pi_2(M) \cong \pi_2(N)$ provided $\Z G$ satisfies
torsion-free cancellation (TFC), which Swan established for certain cyclic and dihedral groups \cite{Sw88}.

The K3 problem list asks for a homeomorphism classification of closed $4$-manifolds with finite fundamental group, and identifies groups of order at most eight as a first goal \cite[Problem~4.50]{BKR26}. One major difficulty is to determine the possible modules $\pi_2(M)$, since this is necessary to begin classifying the equivariant intersection forms and the $k$-invariants \cite[Problems~4.54 and~4.55]{BKR26}. 
The following theorem extends the results obtained for TFC to a substantially larger class, and also incorporates finite groups of order at most eight.

\begin{thmx}\label{thm:finite-all}
Let $G$ be a finite group satisfying at least one of the following conditions:
\begin{clist}{(i)}
\item $G$ has periodic cohomology;
\item $G$ has a cyclic normal subgroup of odd order whose quotient is an ambivalent $2$-group;
\item $G$ has order at most eight.
\end{clist}
If $M$ and $N$ are stably diffeomorphic closed oriented smooth $4$-manifolds with fundamental group $G$, then $\pi_2(M)$ and $\pi_2(N)$ are $\Aut(G)$-isomorphic as $\Z G$-modules. 
\end{thmx}

\begin{remark}
The same conclusion holds for stably homeomorphic closed oriented topological $4$-manifolds: the module extension and stable intersection forms used in the proof are available in this category \cite[Proposition~2.4]{HK88}\cite[Section~2]{HK93}.
\end{remark}

A finite group has \textit{periodic cohomology} if and only if every Sylow subgroup is cyclic or generalised quaternion \cite{Sw60p}. We say that a group is \textit{ambivalent} if every element of $G$ is conjugate to its inverse; examples include all finite Coxeter groups \cite[Theorem 8.7]{Sp74} such as dihedral groups and symmetric groups, and the condition is also preserved by direct products (see \cref{rem:ambivalent}). Condition~\textup{(ii)} includes all ambivalent $2$-groups, by taking the trivial subgroup. It also includes all dihedral groups, by taking the maximum odd order normal subgroup, which is necessarily cyclic.

The proof of \cref{thm:finite-all} applies an algebraic cancellation theorem to Hambleton--Kreck's description of $\pi_2$ by a module extension \cite[Proposition~2.4]{HK88}. Case~\textup{(ii)} also uses stable isometry of the intersection forms.

For $G=C_4\times C_2$ or $C_2^3$, we also prove cancellation of common free summands for arbitrary $\Z G$-lattices, i.e. $P \oplus \Z G \cong Q \oplus \Z G$ implies $P \cong Q$ for all $\Z G$-lattices $P, Q$. We call this property \emph{free summand cancellation} (FSC). Neither group ring has TFC (Theorem~\ref{thm:small-free}), and this gives the first examples of finite groups whose integral group rings have FSC but not TFC.

\subsection{Cancellation for doubles}

The following shows that examples cannot exist among doubles of $2$-complexes whose fundamental groups fall into an even broader range of finite groups.
Given a group presentation, a \textit{commutator twist} is an operation that, given two generators $x, y$ and relators $[x,y]$ and $x^m$, replaces the relator $[x,y]$ by $[x^r,y]$ for some $r$ such that $\gcd(r,m)=1$.  

\begin{thmx}\label{thm:finite-doubles}
Let $X$ and $Y$ be finite $2$-complexes with finite fundamental group $G$ and $\chi(X)=\chi(Y)$. Suppose that at least one of the following holds:
\begin{clist}{(i)}
\item 
$G$ is ambivalent;
\item 
$X$ and $Y$ are homotopy equivalent to presentation $2$-complexes related by a finite sequence of commutator twists.
\end{clist}
Then $\pi_2(M(X))$ and $\pi_2(M(Y))$ are $\Aut(G)$-isomorphic as $\Z G$-modules.
\end{thmx}

A finite abelian group has periodic cohomology precisely when it is cyclic, and is ambivalent precisely when every element has order at most two. For example, $C_p\times C_p$ satisfies neither condition when $p$ is odd. Thus \cref{thm:finite-all} does not cover all finite abelian groups. Nonetheless, combining \cref{thm:finite-doubles} (ii) with results of Metzler \cite{Me76} and Browning \cite[Theorem~1.7]{Bro79} now yields the following.

\begin{corx} \label{cor:abelian}
Let $X$ and $Y$ be finite $2$-complexes with finite abelian fundamental group $G$ and $\chi(X) = \chi(Y)$.
Then $\pi_2(M(X))$ and $\pi_2(M(Y))$ are $\Aut(G)$-isomorphic as $\Z G$-modules.   
\end{corx}

The quadratic bias invariant was introduced by Hambleton--Nicholson as a tool for distinguishing the homotopy types of doubles of $2$-complexes with finite fundamental group. The invariant was used to distinguish doubles whose fundamental groups were finite abelian, as well as a family of the form $Q_8\times(\Z/p)^3$ (see \cite[Section 9]{HN26}).
It follows from \cref{cor:abelian} and another application of \cref{thm:finite-doubles} (ii) that, within each such family, the second homotopy modules are $\Aut(G)$-isomorphic as $\Z G$-modules, where $G$ is the common fundamental group.

\subsection{Non-injectivity of the doubling construction}

It follows from \cref{thm:finite-all} that there exist finite $2$-complexes $X$ and $Y$ with $\chi(X) = \chi(Y)$ and common fundamental group $G$ which are homotopically distinct, such that $\pi_2(M(X))$ and $\pi_2(M(Y))$ are $\Aut(G)$-isomorphic as $\Z G$-modules. 
For example, $Q_{28}$ has periodic cohomology and Mannan--Popiel constructed finite $2$-complexes $X$ and $Y$ with $\chi(X) = \chi(Y)$ with fundamental group $Q_{28}$ which are distinguished by their second homotopy modules \cite[Theorem A]{MP21} (see also \cite{HN27}).

Kreck--Schafer showed that, if $X$ is a finite $2$-complex with $\pi_1(X)$ finite, then
\[ \pi_2(M(X)) \cong \pi_2(X) \oplus \pi_2(X)^* \]
where the involution $\overline g=g^{-1}$ converts the right dual into a left module \cite[Proposition II.4]{KS84}.
Thus the examples of Mannan--Popiel have the property that
\[ \pi_2(X) \not \cong \pi_2(Y) \quad \text{ but } \quad \pi_2(X) \oplus \pi_2(X)^* \cong \pi_2(Y) \oplus \pi_2(Y)^*.\]

This raises the question of whether the doubling construction $X \mapsto M(X)$ is itself injective on homotopy types.
We demonstrate that this is not the case for the examples of Mannan--Popiel.

\begin{thmx}\label{thm:noninjective}
Let $X$ and $Y$ be the presentation $2$-complexes corresponding, respectively, to
\begin{align*}
 \mathcal P_X=\langle x,y\mid y^2x^{-7},\ xyxy^{-1}\rangle,\quad
 \mathcal P_Y=\langle x,y\mid y^2x^{-7},\ 
 y^{-1}xyx^2y^{-1}x^{-2}yx^{-3}
 \rangle.
\end{align*}
Then $X$ and $Y$ have fundamental group $Q_{28}$ and are not homotopy equivalent, but there exists an orientation-preserving homotopy equivalence $M(X)\simeq M(Y)$. In particular, the doubling construction is not injective on homotopy types.
\end{thmx}

The argument does not determine whether these doubles are diffeomorphic. Such a diffeomorphism would also give non-injectivity from generalised Andrews--Curtis classes, or $Q^{**}$-classes, to diffeomorphism classes; compare \cite[Question 1.3]{HN26}.

\subsection{Overview of the proof of \texorpdfstring{\cref{thm:main}}{Theorem A}} \label{ss:proof-idea}

For a left $\Z G$-module $L$, let
\[
 L^*=\Hom_{\Z G}(L,\Z G),\qquad \hmod(L)=L\oplus L^*,
\]
where, as above, the involution $\overline g=g^{-1}$ converts the right dual into a left module. For the groups in this paper, it can be shown that
\[
 \pi_2(M(X))\cong\hmod(\pi_2(X)).
\]
This is false for infinite groups in general, but holds for $G = (\Z/4)^4 * \Z$ (see \cref{lem:double}). The role of cancellation for $L\oplus L^*$ in this construction was discussed by Nicholson \cite[Section 11]{Ni26}.

For \cref{thm:main}, let $H=(\Z/4)^4$ and $G=H*\langle t\rangle$. The proof has four steps.
\begin{clist}{(i)}
\item \emph{Construct the complexes and their stabilisations.}
The presentation complexes $X_n$ differ in a single relator. Generalised Andrews--Curtis moves give
\[
 X_n\vee S^2\simeq_3X_0\vee S^2,
\]
where $\simeq_3$ denotes $3$-deformation equivalence. The spin doubles therefore become diffeomorphic after one $S^2\times S^2$ stabilisation. 

\item \emph{Change coefficients.}
For $n \ge 0$, let $L_n=\pi_2(X_n)$. Abelianisation and localisation at two give
\[
 \rho\colon\Z G\longrightarrow R=\Z_{(2)}[H\times\Z]=S[H],
 \qquad S=\Z_{(2)}[t,t^{-1}].
\]
Let $P_n=R\otimes_{\Z G}L_n$. By Lemma~\ref{lem:base-change},
$R\otimes_{\Z G}\hmod(L_n)\cong\hmod(P_n)$, so it suffices to
distinguish these doubled modules. The same lemma identifies $P_n$
with the kernel of the Fox matrix for $X_n$ over $R$. This matrix is
obtained from the one for $X_0$ by multiplying the last column by
$u_n=2t^n-1$.

\item \emph{Distinguish the doubled modules.}
The column scaling from (ii) induces an isomorphism between
$\hmod(P_n)$ and $\hmod(P_0)$ after extending scalars to
$S[u_n^{-1},\bar u_n^{-1}][H]$ (Lemma~\ref{lem:comparison}).
Thus, any $R$-isomorphism $F \colon \hmod(P_m) \to \hmod(P_n)$, $n \neq m$, induces
an automorphism of $\hmod(P_0)$ after extending scalars to $S[u_n^{-1},\bar u_n^{-1},u_m^{-1},\bar u_m^{-1}][H]$.
To investigate such automorphisms, we pass to $H$-fixed points because
their reduction modulo $4$ identifies with Tate cohomology in degree zero, where we can control the determinant of the induced automorphisms.
The automorphism of the $H$-fixed points of the localisation of $\hmod(P_0)$ induced by $F$ will have determinant $d v_n/v_m$, where $d \in S^\times$ and $v_r=u_r\bar u_r^{-1}$.
By Lemma~\ref{lem:determinant} and the proof of Proposition~\ref{prop:marked}, this determinant must be a monomial modulo $4$. The Laurent coefficients of $dv_n/v_m$ rule this out, so no such $F$ exists. Note that this use of determinants is related to Metzler's bias invariant (see \cite{Me76,La93}).

\item \emph{Allowing for automorphisms of $G$.}
Every automorphism of $G$ descends to one of the finitely many automorphisms of $H\times\Z$. Hence infinitely many of the module classes remain distinct, even modulo the action of $\Aut(G)$.
 \end{clist}

For \cref{thm:finite-all}, Hambleton--Kreck express $\pi_2$, up to a free summand, as an extension of the dual of a third syzygy by the syzygy itself \cite[Proposition~2.4]{HK88}. For periodic groups this extension determines the module up to free stabilisation, and we prove cancellation in its stable class. In case~\textup{(ii)} we also use the stable isometry of the intersection forms at two. For \cref{thm:finite-doubles}, local comparison maps have reduced norm class $u$, and adding the dual changes this to $u\bar u^{-1}$. Ambivalence or the explicit Fox matrix of a commutator twist makes this doubled class trivial. The normal forms for finite abelian groups then give \cref{cor:abelian}.

For \cref{thm:noninjective}, we also compare the intersection forms and $k$-invariants. A transfer argument gives hyperbolic descriptions of the forms. Strong approximation then supplies an integral isometry with enough congruence control to preserve the $k$-invariants. The homotopy equivalence is then obtained from Hambleton--Kreck's result that the quadratic $2$-type determines the homotopy type over groups with $4$-periodic cohomology \cite[Theorem A]{HK88}.

We recall the common geometric construction in \cref{sec:doubles}, prove \cref{thm:main} in \cref{sec:presentations,sec:determinant}, and prove \cref{thm:finite-all,thm:finite-doubles} and \cref{cor:abelian} in \cref{sec:finite}. The example in \cref{thm:noninjective} is treated in \cref{sec:noninjective}. 
From now on, we assume all CW-complexes and manifolds are connected.

\subsection*{Acknowledgements}
TH gratefully acknowledges support by the Deutsche Forschungsgemeinschaft -- Project-ID 286237555 -- TRR 195; and Project-ID 539387714.
JN was supported by the Max Planck Institute for Mathematics.
He would also like to thank F. E. A. Johnson for helpful conversations, and for suggesting the investigation of second homotopy modules of $4$-manifolds.

\subsection*{Declaration of AI usage and origin of the work}

The case of $4$-periodic groups in \cref{thm:finite-all} dates back to the PhD work of the second named author. 
The authors' subsequent collaboration made further progress on this theorem, as part of a long-running project which began in June 2021. 
The original aim of the project was to find closed $4$-manifolds with finite fundamental group and distinct second homotopy modules within a stable diffeomorphism class. The search produced no examples, but laid the foundations for \Cref{thm:finite-all,thm:finite-doubles}, and \cref{cor:abelian}.

The investigation leading to \cref{thm:main} explored many approaches with ChatGPT (GPT-6 Astra), initially with the aim of finding distinct second homotopy modules. Separate attempts to construct an infinite homotopy stable class had repeatedly been unsuccessful. After finding examples with distinct second homotopy modules, the authors examined the construction by hand and found a modification yielding the infinite family. ChatGPT was also used in completing the final details in the proof of \cref{thm:noninjective}, and to aid in creating an initial draft of this manuscript.
The current version has been extensively rewritten and is the responsibility of the authors.

\section{Doubles and second homotopy modules}\label{sec:doubles}

For a finite $2$-complex $X$, choose a smooth $5$-dimensional thickening with trivial stable tangent bundle and write $M(X)$ for its boundary. For presentation complexes, we use handle thickenings built directly from the specified generators and relators. We use models obtained by doubling a $4$-dimensional handle thickening $V\simeq X$:
\[
 M(X)=V\cup_{\partial V}(-V)=\partial(V\times I).
\]
The $2$-handle framings can be chosen so that $V$ is spin. Since $V$ has a $2$-dimensional spine, it is then stably parallelisable. These are the boundary-of-thickening models of Kreck--Schafer \cite[Section II]{KS84}; see also Nicholson \cite[Section 11.1]{Ni26}. In particular, $M(X)$ is closed, oriented, smooth and spin.

\Needspace{10\baselineskip}
\begin{lemma}\label{lem:double}
Let $X$ be a finite $2$-complex and $G=\pi_1(X)$. Then
\[
 \pi_1(M(X))\cong G,\qquad \chi(M(X))=2\chi(X),
\]
and there is an isomorphism of $\Z G$-modules
\[
 \pi_2(M(X))\cong H_2(X;\Z G)\oplus H^2(X;\Z G).
\]
If $H^2(G;\Z G)=H^3(G;\Z G)=0$, it follows that $\pi_2(M(X))\cong\hmod(\pi_2(X))$.
\end{lemma}
\begin{proof}
The fundamental-group and module assertions are \cite[Lemma 11.4]{Ni26}. We use the splitting induced by the map $M(X)\to V\simeq X$. Specifically, its induced map on $\pi_2$ is projection onto $H_2(X;\Z G)$, and its kernel is identified with $H^2(X;\Z G)$ by Poincar\'e--Lefschetz duality. We have $\chi(M(X))=2\chi(V)=2\chi(X)$, since $\chi(\partial V)=0$.

The universal coefficient sequence \cite[proof of Lemma 11.5]{Ni26} is
\[ 0\longrightarrow H^2(G;\Z G)\longrightarrow H^2(X;\Z G)
 \longrightarrow\pi_2(X)^*\longrightarrow H^3(G;\Z G)\longrightarrow0,
\]
where cohomology modules are converted to left modules using the involution. The vanishing assumption gives the final assertion.
\end{proof}

The vanishing in \cref{lem:double} holds for finite groups and for $G=H*\Z$ with $H$ finite:
\begin{equation}\label{eq:groupcohom}
 H^q(G;\Z G)=0 \quad \text{for $q \ge 2$.}
\end{equation}
For finite groups, $\Z G$ is coinduced from the trivial subgroup, so Shapiro's lemma applies. For $H*\Z$, use the cohomology sequence of the Bass--Serre tree \cite[Chapter VII]{Br82}: the finite vertex group has free restricted coefficients, and the other vertex group has cohomological dimension one.

Write $X\simeq_3Y$ if $X$ and $Y$ are related by a $3$-deformation: a sequence of elementary expansions and collapses through complexes of dimension at most three. For presentation complexes, these equivalences are realised by generalised Andrews--Curtis moves \cite[pp.~7--8]{Me76}; see also \cite[Section 2.3]{KMN22}.

\begin{prop}\label{prop:geometric-stabilisation}
Let $X$ and $Y$ be finite $2$-complexes. If $X\vee S^2\simeq_3Y\vee S^2$, then
\[
 M(X)\#(S^2\times S^2)\cong M(Y)\#(S^2\times S^2)
\]
by an orientation-preserving diffeomorphism.
\end{prop}

\begin{proof}
The $3$-deformation equivalence gives generalised Andrews--Curtis moves on presentation models. These preserve the spin boundary of a $5$-dimensional thickening \cite[Proposition~1.5(1),(3)]{Ma11}. Wedging on $S^2$ replaces this boundary by its connected sum with $S^2\times S^2$ \cite[Section~1.2]{KMN22}. This gives the required diffeomorphism between $M(X)\#(S^2\times S^2)$ and $M(Y)\#(S^2\times S^2)$. If this diffeomorphism is not orientation-preserving, then we can obtain one by composing with the reflection exchanging the two halves of the double.
\end{proof}

\section{An infinite family of presentation complexes}\label{sec:presentations}

Let $H=(\Z/4)^4=\langle x,y,z,w\rangle$ and $G=H*\langle t\rangle$, and use the convention $[a,b]=aba^{-1}b^{-1}$. Consider the ordered list
\[
 \mathcal R=\{ x^4y^{-4},\ [x,y]x^4,\ z^4,\ w^4,\
 [x,z],\ [x,w],\ [y,z],\ [y,w],\ [z,w]\}.
\]

Define
\begin{equation}\label{eq:presentations-body}
 \mathcal P_n=\langle x,y,z,w,t\mid\mathcal R,[t^n,x^8]x^4\rangle
 \quad \text{for $n \ge 0$.}
\end{equation}

Let $X_n$ be their presentation complexes, with $L_n=\pi_2(X_n)$. The last relator is
\[
 [t^n,x^8]x^4=t^nx^8t^{-n}x^{-4},
\]
which is $x^4$ when $n=0$. In that case, multiplying the first two relators by the last relator or its inverse, inverting the first, and reordering gives the standard presentation
\[
 \langle x,y,z,w,t\mid x^4,y^4,z^4,w^4,
 [x,y],[x,z],[x,w],[y,z],[y,w],[z,w]\rangle.
\]

\begin{lemma}\label{lem:groups}
Every presentation in \eqref{eq:presentations-body} presents $G$ on the same generating set. For all $n\ge0$, we have
\[
 X_n\vee S^2\simeq_3 X_0\vee S^2.
\]
\end{lemma}

\begin{proof}
Let $c:=x^4$ and $q_n:=t^nc^2t^{-n}c^{-1}$. The common relations give $y^4=c$ and $yxy^{-1}=x^5$. Conjugating the first equality by $y$ gives $x^{20}=x^4$, so $c^4=1$. The relation $q_n=1$ now gives
\[
c=t^nc^2t^{-n}=t^{2n}c^4t^{-2n}=1.
\]
The remaining relations give precisely $H*\langle t\rangle$, and all the displayed relators hold in this group. 

We use elementary relator operations as in \cite[pp.~7--8]{Me76}. Multiplying one relator by a conjugate of another, or its inverse, is a composite of Andrews--Curtis moves. The initial relator is now $q_0=c$.
For each $n$, the preceding group calculation shows that $c$ lies in the normal closure of $\mathcal R$ and $q_n$, while $q_n$ is a product of conjugates of $c$. With one empty relator, this gives
\[
 (\mathcal R,c,1)\sim_Q(\mathcal R,c,q_n)
 \sim_Q(\mathcal R,q_n,c)\sim_Q(\mathcal R,q_n,1),
\]
where $\sim_Q$ denotes equivalence by relator operations with the generators fixed. The first step constructs $q_n$ in the empty slot using $c$; the middle step interchanges the last two slots; the last step erases $c$ using an expression for it as a product of conjugates of the retained relators. Since $q_0=c$, this proves the assertion, using only one additional relator and fixing the generators.
\end{proof}

Let $Y$ be the subcomplex of $X_0$ obtained by omitting the $1$-cell for $t$, so that $X_0 = Y \vee S^1$, and let $J=\pi_2(Y)$. The cellular chain complex $C_*(\wt Y)$ can be augmented to give the exact sequence:
\begin{equation}\label{eq:resolution}
 0\longrightarrow J\longrightarrow \Z H^{10}\xrightarrow{\partial_2^Y} \Z H^4
 \xrightarrow{\partial_1^Y} \Z H\longrightarrow\Z\longrightarrow0.
\end{equation}
Since $X_0 = Y\vee S^1$, it follows from the usual induction formula for homology of a wedge product (see, for example, \cite[Section 6.1]{Ni26}) that:
\begin{equation}\label{eq:induction}
 L_0=\pi_2(X_0)\cong \Z G\otimes_{\Z H}J.
\end{equation}
By \cref{lem:groups}, we have:
\begin{equation}\label{eq:stable}
L_n\oplus \Z G\cong L_0\oplus \Z G \quad \text{for $n \ge 0$.}
\end{equation}

Let $S=\Z_{(2)}[t,t^{-1}]$ and $R=S[H]=\Z_{(2)}[H\times\Z]$, where $\Z_{(2)}$ is the localisation of $\Z$ at the prime two. We equip $R$ with the involution which inverts $t$ and the elements of $H$, and restricts to $t\mapsto t^{-1}$ on $S$. 
Abelianisation and localisation give a map $\rho\colon\Z G\to R$ which respects the involutions on $\Z G$ and $R$.
Let $P_n=R\otimes_{\Z G}L_n$. For $R$-modules, duals are taken over $R$ with the involution defined above, and $\hmod(P)=P\oplus P^*$.

\Needspace{9\baselineskip}
\begin{lemma}\label{lem:base-change}
The modules $P_n$ fit into exact sequences of the form
\begin{equation}\label{eq:Sresolution}
 0\longrightarrow P_n\longrightarrow R^{10}\xrightarrow{\partial_2^{(n)}}R^4
 \xrightarrow{\partial_1} R\longrightarrow S\longrightarrow0.
\end{equation}
Here $\partial_2^{(n)}$ is the second cellular boundary after changing coefficients to $R$, with its zero $t$-coordinate omitted, and $\partial_1$ is independent of $n$.
Furthermore, we have
\[
 P_0\cong S\otimes_{\Z}J,
 \qquad R\otimes_{\Z G}\hmod(L_n)\cong\hmod(P_n).
\]
\end{lemma}
\begin{proof}
Let $U=\im(\partial_2^Y)$ in \eqref{eq:resolution}. All the universal covers of the $X_n$ have the same first boundary $\partial_1^{X}:\Z G^5\to \Z G$, and the images of their second boundaries equal $\ker\partial_1^{X}$.
The cellular complex of $Y\vee S^1$ identifies this image with $\Z G\otimes_{\Z H} U$. Since $R=S\otimes_{\Z}\Z H$, tensoring a $\Z H$-module with $R$ over $\Z H$ is the same as tensoring its underlying abelian group with $S$ over $\Z$. The latter functor is exact: $S$ is a direct sum of copies of the flat $\Z$-module $\Z_{(2)}$. By applying $R \otimes_{\Z H} -$ to the tail of \eqref{eq:resolution}, we obtain an exact sequence:
\begin{equation} \label{eq:res1}
 0\longrightarrow R\otimes_{\Z H}U\longrightarrow R^4
 \xrightarrow{\partial_1} R\longrightarrow S\longrightarrow0.
\end{equation}

Next, from the cellular chain complex $C_*(\wt X_n)$, we obtain an exact sequence:
\[
 0\longrightarrow L_n\longrightarrow \Z G^{10}\longrightarrow \Z G\otimes_{\Z H} U\longrightarrow0.
\]
Since $\Z G$ is free over $\Z H$ and $R$ is flat over $\Z H$, an induced free resolution gives
\[
 \Tor_1^{\Z G}(R,\Z G\otimes_{\Z H} U)\cong\Tor_1^{\Z H}(R,U)=0.
\]
Applying the functor $R \otimes_{\Z G} -$ therefore gives an exact sequence
\begin{equation} \label{eq:res2}
 0\longrightarrow P_n\longrightarrow R^{10}
 \longrightarrow R\otimes_{\Z H}U\longrightarrow0.
\end{equation}

By splicing together \eqref{eq:res1} and \eqref{eq:res2}, we obtain \eqref{eq:Sresolution}: its map out of $R^{10}$ comes from the cellular boundary of $X_n$, while the common tail comes from $Y$. 
Each short exact piece splits over $S$, since its quotient is $S$ tensored over $\Z$ with a free abelian group from \eqref{eq:resolution}.
The formula for $P_0$ follows from \eqref{eq:induction}.

Since $\Z H$ is Noetherian, $J$ is finitely presented. For the dual, the usual finite-presentation argument for Hom and flat base change (compare \cite[Theorem 7.11]{Mat89}) gives, before applying the involution,
\[
 \Hom_{\Z G}(\Z G\otimes_{\Z H}J,\Z G)
 \cong\Hom_{\Z H}(J,\Z G)
 \cong\Hom_{\Z H}(J,\Z H)\otimes_{\Z H} \Z G.
\]
The analogous formula with $R$ identifies the natural map
$
 R\otimes_{\Z G}L_0^*\longrightarrow P_0^*
$
as an isomorphism. The comparison is also an isomorphism for a free module and is compatible with finite direct sums and direct summands. By \eqref{eq:stable}, each $L_n$ is a direct summand of $L_0\oplus\Z G$, so its comparison is an isomorphism too. Hence $R\otimes_{\Z G}\hmod(L_n)\cong\hmod(P_n)$, as required.
\end{proof}

\section{Distinguishing the doubled modules}\label{sec:determinant}

Let $N_H=\sum_{h\in H}h$ be the norm element and let $u_n=2t^n-1$ for $n\ge0$, so $u_0=1$. Order the $2$-cells as the nine relators in $\mathcal R$, followed by $[t^n,x^8]x^4$, in \eqref{eq:presentations-body}; let $e_1,\ldots,e_{10}$ be the corresponding cellular basis. Fox differentiation \cite{Fox53} gives
\begin{equation}\label{eq:Fox}
 \partial_2^{(n)}=\partial_2^{(0)}E_n,\qquad E_n=\operatorname{diag}(1,\ldots,1,u_n).
\end{equation}
Indeed, the last column of $\partial_2^{(0)}$ is $(1+x+x^2+x^3)e_x$, the last column of $\partial_2^{(n)}$ is $u_n$ times this column, and all other columns agree. The $t$-derivatives of all the relators evaluate to zero. By \eqref{eq:Fox}, if $p\in P_n=\ker\partial_2^{(n)}$, then $\partial_2^{(0)}E_np=0$. Hence $E_n$ restricts to a map $P_n\to P_0$.

\subsection{Fixed points and duals}

The comparison maps $E_n$ become isomorphisms after localisation. To detect whether an isomorphism exists before localisation, we restrict to $H$-fixed points, where the modules are free over the commutative ring $S$ and determinants are available. We will compute the comparison determinant and then show that automorphisms of the reference double can change it only by a sixth power modulo four; this is the obstruction described in \cref{ss:proof-idea}.

For an $R=S[H]$-module $P$, let
\[
 P^H=\{v\in P:hv=v\text{ for every }h\in H\},
\]
the $S$-module of $H$-fixed points. If $f$ is an $R$-linear map, write $f^H$ for its restriction to these fixed-point modules.
For an $S$-module $M$, let $\overline{\Hom}_S(M,S)$ denote
the conjugate dual of $M$, consisting of additive maps
$f\colon M\to S$ satisfying $f(bm)=\bar b f(m)$ for all
$b\in S$ and $m\in M$. Its $S$-action is $(bf)(m)=b f(m)$.

The coefficient-of-identity pairing gives the standard identification between the group-ring dual and the conjugate dual over the coefficient ring; compare Curtis--Reiner \cite[Sections 9A and 10D]{CR81}. We record the maps to specify the involution convention. Let $\varepsilon_1:R\to S$ take the coefficient of $1$. For an $R$-module $P$, there is a natural $R$-module isomorphism
\[
 \Psi_P:P^* \longrightarrow \overline{\Hom}_S(P,S),
 \quad F\longmapsto\bigl(m\longmapsto\overline{\varepsilon_1(F(m))}\bigr),
\]
where $H$ acts on the conjugate dual by $(hf)(m)=f(h^{-1}m)$. Its inverse sends $f$ to the $R$-linear functional
\[
 m\longmapsto\sum_{h\in H}\overline{f(h^{-1}m)}\,h.
\]

By abuse of notation, we write $\Psi$ for $\Psi_P$ for every $P$.

\begin{lemma}\label{lem:fixed}
The following hold for all $n \ge 0$:
\begin{clist}{(i)}
  \item
    Each $P_n^H$ is free of rank six over $S$, with basis $N_H e_5,\ldots,N_H e_9,\quad N_H(e_{10}-u_n e_2)$. 
   \item
     In the bases in \textup{(i)}, the map $E_n^H:P_n^H\to P_0^H$ has matrix $\operatorname{diag}(1,1,1,1,1,u_n)$.
    \item We have $N_H P_n=4P_n^H$ and $N_H P_n^*=4(P_n^*)^H$.
    \item
Restriction after $\Psi$ induces an $S$-module isomorphism
\[
 (P_n^*)^H \longrightarrow 64\,\overline{\Hom}_S(P_n^H,S),
 \quad F\longmapsto\Psi(F)|_{P_n^H}.
\]
Writing $p_{n,1},\ldots,p_{n,6}$ for the basis in \textup{(i)}, the module $(P_n^*)^H$ has basis $F_{n,1},\ldots,F_{n,6}$ characterised by
\[
 \Psi(F_{n,i})(p_{n,j})=64\delta_{ij}.
\]
\end{clist}
\end{lemma}
\begin{proof}
(i) and (ii).
Let $\eps_H:R=S[H]\to S$ be the augmentation.
An element of $R$ is $H$-fixed precisely when its coefficients in the $H$-basis are all equal, so $R^H=SN_H$. Moreover, $hN_H=N_H$ for every $h\in H$, and hence
$
 aN_H=\eps_H(a)N_H
$
for $a \in R$.
Thus, in the bases $N_H e_1,\ldots,N_H e_{10}$ of $(R^{10})^H$ and $N_H e_x,N_H e_y,N_H e_z,N_H e_w$ of $(R^4)^H$, the matrix of $(\partial_2^{(n)})^H$ is obtained by applying $\eps_H$ to each entry of $\partial_2^{(n)}$, yielding
\[
\begin{pmatrix}
 4&4&0&0&0&0&0&0&0&4u_n\\
 -4&0&0&0&0&0&0&0&0&0\\
 0&0&4&0&0&0&0&0&0&0\\
 0&0&0&4&0&0&0&0&0&0
\end{pmatrix}.
\]
Since $P_n^H=\ker((\partial_2^{(n)})^H)$, the kernel equations give the basis in (i). The map $E_n$ fixes $N_H e_5,\ldots,N_H e_9$ and sends the last basis vector to $u_nN_H(e_{10}-e_2)$, proving (ii).

(iii) and (iv).
Write $\HH^q$ for Tate cohomology and consider $R = S[H]$-modules as $\Z[H]$-modules. Here $H$ acts trivially on $S$ via $\eps_H$.
Since $R\cong\Z[H]\otimes_{\Z}S$ is induced from the trivial
subgroup, the three middle terms of \eqref{eq:Sresolution}
have zero Tate cohomology.
The connecting morphisms of the long exact Tate sequence applied to \eqref{eq:Sresolution} yield
\[
 P_n^H/N_H P_n=\HH^0(H;P_n)
 \cong\HH^{-3}(H;S)=H_2(H;S)\cong(S/4)^6.
\]
The last isomorphism uses flatness of $S$ over $\Z$ and the K\"unneth formula: $H_1(C_4;\Z)\cong\Z/4$ and $H_2(C_4;\Z)=0$, so each of the six pairs of cyclic factors contributes one copy of $\Z/4$ to $H_2(H;\Z)$; see \cite[Section V.6]{Br82}. Thus $4P_n^H\subseteq N_H P_n$. The induced surjection $P_n^H/4P_n^H\to P_n^H/N_H P_n$ is between free modules of the same finite rank over the commutative ring $S/4S$, so it is an isomorphism. This proves the first equality in (iii).

For $F\in(P_n^*)^H$, write $f=\Psi(F)$. Since $f$ is $H$-invariant, $f(N_H m)=256f(m)$. Together with $N_H P_n=4P_n^H$, this gives $f(P_n^H)\subseteq64S$. Conversely, for $g\in\overline{\Hom}_S(P_n^H,S)$, the element
\[
 F=\Psi^{-1}\bigl(m\longmapsto g(N_H m/4)\bigr)
\]
lies in $(P_n^*)^H$ and satisfies $\Psi(F)|_{P_n^H}=64g$. The map in (iv) is injective because $\Psi(F)|_{P_n^H}=0$ implies $256\Psi(F)=0$, hence $F=0$ by torsion-freeness of $S$. The $F_{n,i}$ correspond to the scaled coordinate functionals under this isomorphism, proving the basis assertion in (iv).

Under the preceding identification of the $R$-dual with the conjugate $S$-dual, dualising the $S$-split sequence \eqref{eq:Sresolution} gives an exact sequence
\[
 0\longrightarrow S\longrightarrow R\longrightarrow R^4
 \longrightarrow R^{10}\longrightarrow P_n^*\longrightarrow0.
\]
The three middle terms again have zero Tate cohomology. The connecting isomorphisms now shift degree zero at the right-hand end to degree three at the left-hand end, giving
\[
 \HH^0(H;P_n^*)\cong\HH^3(H;S)=H^3(H;S)\cong(S/4)^6.
\]
For the last equality, the universal coefficient theorem \cite[Section III.1, Exercise 3]{Br82} gives
\[
 H^3(H;S)\cong\Ext^1_{\Z}(H_2(H;\Z),S),
\]
since $H_3(H;\Z)$ is finite and $S$ is torsion-free.
Since $H_2(H;\Z)\cong(\Z/4)^6$ and
$\Ext^1_{\Z}(\Z/4,S)\cong S/4S$, we obtain
$\Ext^1_{\Z}(H_2(H;\Z),S)\cong(S/4S)^6$.
By (iv), $(P_n^*)^H$ is free of rank six over $S$, so the same
equal-rank argument proves the second equality in (iii).
\end{proof}

Recall that $S=\Z_{(2)}[t,t^{-1}]$ and $R=S[H]$.
Fix a finite set $I\subseteq\Z_{\ge0}$. We now localise $S$ so that the maps $E_n^H$ become isomorphisms for $n\in I$.
It suffices to invert the corresponding $u_n$. Put
\[
 C=S[u_n^{-1},\overline{u_n}^{-1}:n\in I].
\]
We also invert the $\overline{u_n}$ so that the set of denominators is stable under the involution on $S$. Thus the involution extends to $C$ by $\overline{a/b}=\bar a/\bar b$, and to $C[H]$ by $\bar h=h^{-1}$ for $h\in H$.

Recall that $J=\pi_2(Y)$, where $Y$ is obtained from $X_0$ by deleting the free circle $t$.
Let $J_C=C\otimes_{\Z}J$ and $P_{n,C}=C\otimes_S P_n$.

\begin{lemma}\label{lem:comparison}
The following hold for all $n\in I$:
\begin{clist}{(i)}
\item The map $E_n$ induces a $C[H]$-module isomorphism $P_{n,C}\to J_C$, and hence an isomorphism
\[
  \mathcal{E}_n=E_n\oplus(E_n^{-1})^*:
 \hmod(P_{n,C}) \longrightarrow \hmod(J_C).
\]
\item In the localised bases from \cref{lem:fixed}\textup{(i)} and the corresponding scaled dual bases from part~\textup{(iv)} of that lemma, we have
\[
 \det(\mathcal{E}_n^H)=v_n:=u_n\bar u_n^{-1},\qquad v_0=1.
\]
\end{clist}
\end{lemma}
\begin{proof}
Since $C$ is flat over $S$, tensoring with $C$ preserves kernels.
For $r=0,n$, let $D_r=C\otimes_S\partial_2^{(r)}$.
Then $P_{r,C}=\ker D_r$.

Since $H$ is finite, $P_r^H$ is the kernel of
$P_r\to\bigoplus_{h\in H}P_r$, $m\mapsto((h-1)m)_h$.
Moreover, $P_r\subseteq R^{10}$ is finitely presented because $R$ is Noetherian
\cite[Theorem~3.1(iii)]{Mat89}.
Dualising a finite presentation $R^a\to R^b\to P_r\to0$
similarly expresses $P_r^*$ as a kernel between finite free modules.
The usual Hom comparison under flat base change applies to $R\to C[H]$
\cite[Theorem~7.11]{Mat89}.
As the involution extends to $C[H]$, we therefore obtain
\[
 C\otimes_S P_r^H\cong(P_{r,C})^H,\qquad
 C\otimes_S P_r^*\cong(P_{r,C})^*.
\]

(i)
Since $u_n\in C^\times$, the matrix $E_n$ is invertible over $C$.
Equation \eqref{eq:Fox} gives $D_n=D_0E_n$, so $E_n$ restricts
to an isomorphism
$P_{n,C}=\ker D_n\longrightarrow\ker D_0=P_{0,C}=J_C$.

(ii)
By \cref{lem:fixed}(ii), with respect to the chosen basis, the matrix corresponding to $E_n^H$ is equal to
\[
 A_n:=\operatorname{diag}(1,1,1,1,1,u_n).
\]
For $((E_n^{-1})^*)^H \colon (P_{n,C}^\ast)^H \to (J_C^\ast)^H$, with respect to the dual basis of \cref{lem:fixed}(iv), the matrix is given by
\[
 \overline{A_n^{-1}}^{\,T}
 =\operatorname{diag}(1,1,1,1,1,\bar u_n^{-1}).
\]
Note that the normalisation $\Psi(F_{n,i})(p_{n,j})=64\delta_{ij}$ is the same for domain and codomain, so the common factor $64$ cancels when computing the second matrix. Multiplying the determinants proves (ii).
\end{proof}

\subsection{Arbitrary endomorphisms of the doubled module}

We now use the comparison maps of \cref{lem:comparison} to
distinguish the doubled modules over $R$.
For $I=\{m,n\}$, an $R$-isomorphism
$\hmod(P_m)\to\hmod(P_n)$ would yield an automorphism of
$\hmod(J_C)$.
We first show that the determinant of any $C[H]$-automorphism
of $\hmod(J_C)$ on $H$-fixed points is a sixth power of a unit
modulo four, and then use the comparison determinants to rule
out such an isomorphism when $m\ne n$.

\begin{lemma}\label{lem:determinant}
Let $\theta$ be a $C[H]$-endomorphism of $J_C\oplus J_C^*$.
The following hold:
\begin{enumerate}[\normalfont (i)]
\item There exist $a,b\in C/4$ such that, for any choice of
$C/4$-bases of $\HH^0(H;J_C)$ and $\HH^0(H;J_C^*)$,
the endomorphism induced by $\theta$ on their direct sum has matrix
\[
 \begin{pmatrix}aI_6&\ast\\0&bI_6\end{pmatrix} \in M_{12}(C/4).
\]
\item With $a,b$ as in \textup{(i)}, we have
\[
 \det(\theta^H)\equiv(ab)^6\bmod4.
\]
\item If $\theta$ is invertible, then $a,b\in(C/4)^\times$.
\end{enumerate}
\end{lemma}
\begin{proof}
(i)
Tensoring \eqref{eq:resolution} with $C$ gives the exact sequence
\[
 0\longrightarrow J_C\longrightarrow C[H]^{10}
 \longrightarrow C[H]^4\longrightarrow C[H]
 \longrightarrow C\longrightarrow0,
\]
whose short exact pieces split over $C$.
Finite free $C[H]$-modules are injective relative to $C$-split inclusions \cite[Section VI.2]{Br82}. Thus every endomorphism of $J_C$ extends successively to an endomorphism of this exact sequence, inducing multiplication by a scalar on its final term $C$. Naturality of dimension shifting gives the same scalar action on
$\HH^0(H;J_C)\cong\HH^{-3}(H;C)=H_2(H;C)$.

Under $\Psi$, dualising the sequence is the same as taking its conjugate $C$-dual, so the dual sequence is also exact. Its middle terms are free over $C[H]$. Their projectivity lets us lift any endomorphism of $J_C^*$ successively through the dual sequence, inducing a scalar endomorphism of its initial term $C$. Dimension shifting then gives scalar action on $\HH^0(H;J_C^*)\cong\HH^3(H;C)=H^3(H;C)$.

We now show that every map $f:J_C\to J_C^*$ induces zero on $\HH^0$. For each coordinate pair $i<j$, put
\[
 H_{ij}=\langle x_i,x_j\rangle\cong(\Z/4)^2,
 \qquad (x_1,x_2,x_3,x_4)=(x,y,z,w).
\]
Under the natural identification
\[
 H^3(H;C)\cong\Ext^1_{\Z}\bigl(H_2(H;\Z),C\bigr)
 \cong(C/4)^6,
\]
the K\"unneth decomposition above identifies the inclusions of the six $H_{ij}$ with the six homology summands. By naturality of the universal coefficient theorem, restriction in cohomology therefore gives the six coordinate projections
\cite[Section III.1, Exercise 3, and Theorem V.6.4(iii)]{Br82}.
Thus these restrictions are jointly injective.

Fix $K=H_{ij}$. The resolution and its dual remain exact with free middle terms on restriction to $K$. Naturality of dimension shifting gives a commutative square
\[
\begin{CD}
 H_2(H;C) @>{\HH^0(H;f)}>> H^3(H;C)\\
 @V{\operatorname{tr}_K^H}VV @VV{\res_K^H}V\\
 H_2(K;C) @>{\HH^0(K;\res f)}>> H^3(K;C).
\end{CD}
\]
Both vertical maps come from Tate restriction: in degree $-3$ it is homological transfer, and in degree $3$ it is ordinary cohomological restriction \cite[Sections III.9 and VI.5]{Br82}.

Let $\iota:K\hookrightarrow H$ be the inclusion and $p:H\to K$ the coordinate projection. Since $p\iota=\id_K$ and $\iota_*\operatorname{tr}_K^H=[H:K]\id$ \cite[Proposition~III.9.5(ii)]{Br82}, we have
\[
 \operatorname{tr}_K^H
 =p_*\iota_*\operatorname{tr}_K^H
 =16p_*=0,
\]
because $H_2(K;C)\cong C/4$. The square therefore gives
$\res_K^H\,\HH^0(H;f)=0$ for every $H_{ij}$, which implies $\HH^0(H;f)=0$.

(ii)
Let $M=J_C\oplus J_C^*$. By \cref{lem:fixed}(iii), after localising from $S$ to $C$ (which yields a flat extension $S \to C$),
$\HH^0(H;M)=M^H/N_H M=M^H/4M^H$.
Hence the map in (i) is $\theta^H$ reduced modulo 4, so taking determinants gives $\det(\theta^H)\equiv(ab)^6\bmod4$.

(iii)
If $\theta$ is invertible, $(ab)^6$ is a unit, so both $a$ and $b$ are units.
\end{proof}

\begin{prop}\label{prop:marked}
The modules $\hmod(P_n)$, $n\ge0$, are pairwise non-isomorphic over $R$. Consequently, the $\Z G$-modules $\hmod(L_n)$ are pairwise non-isomorphic, even after localising at two.
\end{prop}
\begin{proof}
  Suppose $F \colon \hmod(P_m)\to\hmod(P_n)$ is an isomorphism with $0\le m<n$. By \cref{lem:fixed}, $\det F^H\in S^\times = \Z_{(2)}[t, t^{-1}]^\times $, so $\det F^H=ct^j$ for some $c\in\Z_{(2)}^\times$ and $j\in\Z$. Take $I=\{m,n\}$ in the definition of $C$. Applying \cref{lem:comparison,lem:determinant} to the automorphism $\mathcal{E}_n F_C\mathcal{E}_m^{-1}$ of $\hmod(J_C)$ gives
\[
 ct^j\frac{v_n}{v_m}\equiv w^6\bmod4
 \qquad\text{for some }w\in(C/4)^\times.
\]

Every $u_r$ and $\overline{u_r}$ is already a unit modulo 4, so
\[
 C/4\cong(\Z/4)[t,t^{-1}].
\]
Reducing $w$ modulo two gives a unit of $\mathbb F_2[t,t^{-1}]$, hence $t^l$ for some $l\in\Z$, since the difference between the highest and lowest exponents is nonnegative and additive under multiplication.
Thus $w=t^l(1+2f(t))$ for some $f\in(\Z/4)[t,t^{-1}]$.
As $(1+2f(t))^2=1+4f(t)+4f(t)^2=1$ modulo 4, $w^6=t^{6l}$.
Thus $v_n/v_m$ must be a signed monomial modulo 4. Using $(2a-1)^{-1}\equiv2a-1\bmod4$ and discarding terms divisible by 4 gives
\[
\begin{aligned}
 \frac{v_n}{v_m}
 &=\frac{(2t^n-1)(2t^{-m}-1)}{(2t^{-n}-1)(2t^m-1)}\\
 &\equiv(2t^n-1)(2t^{-m}-1)(2t^{-n}-1)(2t^m-1)\\
 &\equiv1+2(t^n+t^{-n}+t^m+t^{-m})\bmod4.
\end{aligned}
\]
Its coefficient at $t^n$ is $2$, which is impossible for a signed monomial.

The last assertion follows from \cref{lem:base-change}, since $\rho\colon\Z G\to R$ factors through $\Z_{(2)}G$.
\end{proof}

\subsection{\texorpdfstring{Allowing for automorphisms of $G$}{Allowing automorphisms of G}}

We write $U\cong_{\Aut(G)}V$, and call $U$ and $V$ $\Aut(G)$-isomorphic, if the two $\Z G$-modules become isomorphic after twisting one action by an automorphism of $G$.
Equivalently, there exist $\alpha\in\Aut(G)$ and an $\alpha$-semilinear
isomorphism $F\colon U\to V$, that is, a $\Z$-linear bijection satisfying
$F(gx)=\alpha(g)F(x)$ for all $g\in G$ and $x\in U$.

\begin{prop}\label{prop:unmarked}
The family $\{\hmod(L_n):n\ge0\}$ contains infinitely many $\Aut(G)$-isomorphism classes.
\end{prop}
\begin{proof}
Every $\alpha\in\Aut(G)$ induces an automorphism of
$G_{\mathrm{ab}}=H\times\Z$, and hence an automorphism
$\widetilde\alpha$ of $R$. Thus $\widetilde\alpha\circ\rho=\rho\circ\alpha$ on $\Z G$. By \cite[Theorem~1]{Al62}, $\Aut(H\times\Z)$ is finite, since $\langle t\rangle$ is a central cyclic subgroup of finite index.

An $\alpha$-semilinear isomorphism $F:\hmod(L_m)\to\hmod(L_n)$ induces a $\widetilde\alpha$-semilinear isomorphism
\[
R\otimes_{\Z G}\hmod(L_m)\longrightarrow R\otimes_{\Z G}\hmod(L_n),
 \qquad r\otimes x\longmapsto
 \widetilde\alpha(r)\otimes F(x).
\]
By \cref{lem:base-change}, its source and target are
$\hmod(P_m)$ and $\hmod(P_n)$. Thus each
$\Aut(G)$-isomorphism class in our family specialises
into a single $\Aut(H\times\Z)$-orbit.

By \cref{prop:marked}, the modules $\hmod(P_n)$ are
pairwise non-isomorphic. Each such orbit therefore
contains the isomorphism class of $\hmod(P_n)$ for at most
$\lvert\Aut(H\times\Z)\rvert$ values of $n$,
so infinitely many $\Aut(G)$-isomorphism classes remain.
\end{proof}

\begin{proof}[Proof of \cref{thm:main}]
By \cref{prop:unmarked}, there is an infinite set $\mathcal I\subseteq\Z_{\ge0}$ such that $\hmod(L_m)\not\cong_{\Aut(G)}\hmod(L_n)$ for distinct $m,n\in\mathcal I$. Take spin doubles $M(X_n)$ for these indices and reindex them as $M_1,M_2,\ldots$. By \eqref{eq:groupcohom} and \cref{lem:double}, they have fundamental group $G$, Euler characteristic $2\chi(X_n)=12$, and the required second homotopy modules. They are therefore pairwise not homotopy equivalent, with either choice of orientation.

For every $n$, \cref{lem:groups} gives $X_n\vee S^2\simeq_3X_0\vee S^2$ by generalised Andrews--Curtis moves. Choose the presentation handle models specified in \cref{sec:doubles}. Proposition~\ref{prop:geometric-stabilisation} gives orientation-preserving diffeomorphisms
$M(X_n)\#(S^2\times S^2)\cong M(X_0)\#(S^2\times S^2)$ for $n\ge0$.
\end{proof}

\section{Finite fundamental groups}\label{sec:finite}

Throughout this section, we let $G$ be a finite group. A $\Z G$-lattice is a finitely generated $\Z G$-module which is free as an abelian group. For such a module $M$, taking the coefficient of the identity identifies the group-ring dual with the integral contragredient dual:
\[
 M^*\cong\Hom_{\Z}(M,\Z),\qquad
 (g\varphi)(m)=\varphi(g^{-1}m).
\]
After scalar extension, we use the contragredient dual
over the corresponding coefficient ring, with the natural
identifications
\[
 (L^*)_{\Q}\cong(L_{\Q})^*,\qquad
 (L^*)_p\cong(L_p)^*.
\]

\subsection{Genus and reduced norms}
We write $L_{\Q}=\Q\otimes_{\Z}L$ and $L_p=\Z_p\otimes_{\Z}L$, where $\Z_p$ denotes the ring of $p$-adic integers. Lattices are in the same \emph{genus} if their completions are isomorphic at every prime. Cancellation holds over $\Z_pG$ by Krull--Schmidt \cite[Theorem 20.6.2]{Vo21}. For a lattice $L$, we put
\[
 E=\End_{\Q G}(L_{\Q}),\qquad K=Z(E),\qquad
 E_p(L)=\End_{\Z_pG}(L_p).
\]
For $F=\Q$ or $\Q_p$, we regard reduced norms as elements
of $Z(FG)^\times$, extending them by $1$ on the simple components
absent from the module. We likewise identify the unit group of
a central direct factor $eZ(FG)$ with its image in $Z(FG)^\times$
under $a\mapsto a+(1-e)$, and make the corresponding identification
componentwise for finite ideles. Under these identifications,
projections onto central factors, including $\pi_L$ and $\pi_{L,p}$
below, retain the indicated components and replace the others
by $1$. Subscripts indicating components of central units use
the same convention.

The centre $K$ is a product of number fields with maximal order $\mathcal{O}_K$. The algebra $E$ satisfies the Eichler condition if no simple component of $E$ is isomorphic to a totally definite quaternion algebra. Under this condition, Jacobinski's genus theorem \cite[Theorem~2.2]{Jac68} (also discussed in \cite[Sections 49--51]{CR87}) identifies the isomorphism classes in the genus of $L$ with
\begin{equation}\label{eq:normclass}
 \mathcal{G}(L)=J_f(K)\Big/\left(\nrd(E^\times)
             \prod_p\nrd(E_p(L)^\times)\right),
\end{equation}
where the class of $L$ corresponds to the identity.
Here, $J_f(K)=\prod_p'(K\otimes_{\Q}\Q_p)^\times$ is the finite idele group, restricted with respect to $(\mathcal O_K\otimes\Z_p)^\times$.
The notation $\nrd$ denotes the reduced norm on $E$ and on
$E\otimes_{\Q}\Q_p$, with the latter restricted to $E_p(L)^\times$
in \eqref{eq:normclass}.
For a $\Z G$-lattice $L'$ in the genus of $L$,
we fix a $\Q G$-isomorphism $L'_\Q\cong L_\Q$ and use it to regard $L'$ as a lattice in $L_\Q$.
If $(\alpha_p)_p$ is a family of isomorphisms $\alpha_p \colon L_p \to L'_p$ of $\Z_pG$-modules, then the class of $L'$ in $\mathcal{G}(L)$ corresponds to the coset of the finite idele $(\nrd(\alpha_p))_p\in J_f(K)$, where each $\alpha_p$ is regarded as an element of $(E\otimes_{\Q}\Q_p)^\times$ by extending scalars to $\Q_p$.
Changing the rational identification multiplies this tuple by an element of $\nrd(E^\times)$, so its class in $\mathcal{G}(L)$ is independent of this choice.

For a $\Z_pG$-lattice $U$, the centre of $\End_{\Q_pG}(\Q_p\otimes_{\Z_p}U)$ is canonically a direct summand of $Z(\Q_p G)$.
We write $\nu_p(U)\subseteq Z(\Q_pG)^\times$ for the image
of $\Aut_{\Z_pG}(U)$ under the reduced norm, using the standing
convention.
Thus all local norm images lie in one fixed group.

\begin{lemma}[Local norm images]\label{lem:local-norm-sums}
For $\Z_pG$-lattices $U,V$ and $r\ge1$, we have
\begin{equation}\label{eq:local-norm-sums}
 \nu_p(U\oplus V)=\nu_p(U)\nu_p(V),
 \qquad \nu_p(U^r)=\nu_p(U).
\end{equation}
In particular, if $V$ is a direct summand of a finite sum of copies of $U$, then $\nu_p(V)\subseteq\nu_p(U)$.
\end{lemma}

\begin{proof}
The rings $\End_{\Z_pG}(U)$ and $\End_{\Z_pG}(V)$
are finitely generated as $\Z_p$-modules and hence semilocal.
Swan's decomposition of automorphisms \cite[Proposition 4.3]{Sw88} therefore reduces their norms to those of block-diagonal and elementary automorphisms. The reduced norms of the block-diagonal automorphisms form the subgroup $\nu_p(U)\nu_p(V)$, while elementary automorphisms have reduced norm $1$. This proves the first equality. The second follows by iteration, since $\nu_p(U)$ is a subgroup. The final assertion follows by extending an automorphism of $V$ by the identity on a complement.
\end{proof}

\begin{prop}[Cancellation by reduced norms]\label{prop:norm-cancellation}
Let $L,T$ be $\Z G$-lattices with $L\ne0$, and suppose that $E=\End_{\Q G}(L_{\Q})$ satisfies the Eichler condition. We let $K=Z(E)$ and write $\pi_L:Z(\Q G)^\times\to K^\times$ for the canonical projection, with completions $\pi_{L,p}:Z(\Q_pG)^\times\to(K\otimes\Q_p)^\times$. Assume that
\[
 \pi_{L,p}\bigl(\nu_p(L_p\oplus T_p)\bigr)
 =\nu_p(L_p)
\]
for all rational primes $p$.
Then, for every $\Z G$-lattice $N$ and every $r\ge0$,
\[
 L\oplus T^r\cong N\oplus T^r
 \quad\Longrightarrow\quad L\cong N.
\]
\end{prop}

\begin{proof}
The case $r=0$ is immediate, so we assume $r\ge1$. Local cancellation shows that $L$ and $N$ lie in the same genus. We fix a rational identification $N_{\Q}\cong L_{\Q}$ and local isomorphisms $a_p:L_p\to N_p$. In these coordinates, we write $f$ for the stable isomorphism and let $\eta=\pi_L(\nrd(f))$.

Over each simple component of $K$, adding $T_{\Q}^r$ only increases matrix sizes over the same division algebras.
The image of the reduced norm is unchanged by increasing these matrix sizes \cite[III.1.2.4]{We13}, so $\eta\in\nrd(E^\times)$.
By \cref{lem:local-norm-sums} and the hypothesis, we have $\pi_{L,p}(\nu_p(L_p\oplus T_p^r))=\nu_p(L_p)$. Since $f^{-1}(a_p\oplus\id_{T_p^r})$ is an automorphism of $L_p\oplus T_p^r$, we obtain
\[
 \nrd(a_p)\eta^{-1}\in\nu_p(L_p)
\]
for every rational prime $p$.
Dividing $(\nrd(a_p))_p$ by the global norm $\eta$ therefore gives a family of local automorphism norms. Hence the class of $N$ in \eqref{eq:normclass} is trivial, and $L\cong N$.
\end{proof}

We now consider the situation of doubled modules.
We will repeatedly use the following consequence of Jacobinski's genus theorem to establish isomorphisms between doubled modules.

\begin{lemma}\label{lem:double-genus}
For $\Z G$-lattices $M,N$ in the same genus,
we fix a $\Q G$-isomorphism $N_\Q\cong M_\Q$ and use it to regard $N$ as a lattice in $M_\Q$.
We use the inverse dual isomorphism to regard $N^*$ as a lattice in $M_\Q^*$.
We choose local isomorphisms $\alpha_p:M_p\to N_p$ and let $u=(\nrd(\alpha_p))_p$.
We let $K=Z(\End_{\Q G}(\hmod(M)_\Q))$ and write $z\mapsto\bar z$ for its involution induced by $g\mapsto g^{-1}$.
Then the class of $\hmod(N)$ in $\mathcal{G}(\hmod(M))$ is represented by $u\bar u^{-1}$.
In particular, $\hmod(M)\cong\hmod(N)$ if and only if $u\bar u^{-1}$ represents the identity in $\mathcal{G}(\hmod(M))$.
\end{lemma}

\begin{proof}
Every finite-dimensional $\Q G$-module is self-dual. Thus the irreducible multiplicities in $\Q\otimes\hmod(M)$ are even, and its endomorphism algebra satisfies the Eichler condition.
For each prime $p$, the map
$\alpha_p\oplus(\alpha_p^{-1})^*:M_p\oplus M_p^*\to N_p\oplus N_p^*$
is a $\Z_pG$-isomorphism. Its scalar extension to $\Q_p$ has
reduced norm $\nrd(\alpha_p)\overline{\nrd(\alpha_p)}^{-1}$,
by compatibility of reduced norms with duality.
The assertion now follows from \eqref{eq:normclass}.
\end{proof}

\subsection{Commutator powers}

\begin{lemma}\label{lem:power}
We take $\mathcal P$ to be a finite presentation of a finite group $G$ in which $x$ and $y$ are two distinct specified generators and $x^m$ and $[x,y]$ are relators, with $m\ge1$. The presentation $\mathcal P_r$ is obtained by replacing the latter relator by $[x^r,y]$, leaving the other relators unchanged, where $(r,m)=1$. Their presentation complexes satisfy
\[
 \hmod(\pi_2(X_{\mathcal P_r}))\cong\hmod(\pi_2(X_{\mathcal P})).
\]
The isomorphism respects the common polarisation by the displayed generators.
\end{lemma}

\begin{proof}
The relation $x^m=1$ implies that commuting with $x^r$ is equivalent to commuting with $x$, so the $G$-polarisation is unchanged. Replacing $r$ by $r+km$ is realised by elementary relator operations using the retained relator $x^m$. The Chinese remainder theorem therefore lets us assume $r>0$ and $(r,|G|)=1$.

We let $L_i=\pi_2(X_{\mathcal P_i})$, $C=\langle x\rangle$, $T=\Z[G/C]$, and $s=1+x+\cdots+x^{r-1}$. In \eqref{eq:Fox} we used columns over the commutative ring $R$. Here we write elements of free left $\Z G$-modules as rows, so their boundary matrices act on the right.

We write $F$ for the free $\Z G$-module on the relators, with the changed relator first.
For $i\in\{1,r\}$, the map $\partial_2^{(i)}:F\to C_1$
is the second cellular boundary map of the universal cover
of $X_{\mathcal P_i}$, where $\mathcal P_1=\mathcal P$
and $C_1$ is the free $\Z G$-module on the generators.
The common first cellular boundary map $\partial_1:C_1\to\Z G$ is defined by $\partial_1(e_g)=g-1$ on the basis elements corresponding to the displayed generators $g$.
Thus $L_i=\ker\partial_2^{(i)}$.
We order the generators with $x,y$ first and write $B_i$ for the matrix of $\partial_2^{(i)}$, so $\partial_2^{(i)}(v)=vB_i$ for $v\in F$.
Using the convention $[a,b]=aba^{-1}b^{-1}$, the first rows of these matrices, obtained by evaluating the Fox derivatives in $\Z G$, are
\[
 b_1=(1-y,x-1,0,\ldots,0),\qquad
 b_r=(s(1-y),x^r-1,0,\ldots,0)=s b_1,
\]
where we use $xy=yx$ in $G$ and $x^r-1=s(x-1)$.
All remaining rows agree, so $B_r=DB_1$ for $D=\operatorname{diag}(s,1,\ldots,1)$. Hence
\begin{equation}\label{eq:commutator-fox}
 \partial_2^{(r)}=\partial_2^{(1)}\circ d,
\end{equation}
where $d:F\to F$ is given by $d(v)=vD$. In $\Z C/(s)$, the identity $(x-1)s=x^r-1$ gives $x^r=1$. Since $(r,|C|)=1$, this forces $x=1$, and then $s=r$. Hence $\Z C/(s)\cong\Z/r\Z$ with trivial $C$-action. Inducing from $C$ to $G$ gives isomorphisms of left $\Z G$-modules
\[
 \Z G/\Z Gs\cong\Z G\otimes_{\Z C}\bigl(\Z C/(s)\bigr)
 \cong\Z G\otimes_{\Z C}(\Z/r\Z)\cong T/rT.
\]
Here we use $\Z G\otimes_{\Z C}\Z\cong\Z[G/C]=T$.
In particular, $\operatorname{coker}(d)\cong T/rT$ is finite. Since $d$ becomes an automorphism after tensoring with $\Q$ and $F$ is torsion-free over $\Z$, the map $d$ is injective. Both Fox maps have image $J=\ker\partial_1$, and $d$ induces the identity on $J$. The snake lemma gives
\[
 0\longrightarrow L_r\xrightarrow{f}L_1
 \longrightarrow T/rT\longrightarrow0,
\]
where $f=d|_{L_r}$.
We use $f$ to identify the rational modules. For $p\nmid r$, the lattices $f((L_r)_p)$ and $(L_1)_p$ are equal.

For $p\mid r$, every $\Z_pG$-lattice is projective since $p\nmid |G|$ \cite[(19.5)(ix)]{CR81}.
Projectivity of $(L_1)_p$ lifts the quotient map
$(L_1)_p\to T_p/rT_p$ to a map $h:(L_1)_p\to T_p$.
Since $h$ is surjective modulo $r$ and $p\mid r$,
Nakayama's lemma shows that $h$ is surjective.
Projectivity of $T_p$ then gives a splitting of $h$.
Thus, in suitable coordinates,
\[
 (L_1)_p=T_p\oplus U_p,\qquad
 f((L_r)_p)=rT_p\oplus U_p.
\]
We choose the $\Z_pG$-isomorphism $a_p=\operatorname{diag}(r,1)$ from the first lattice to the second. The permutation lattice $T_p$ is self-dual, so compatibility of reduced norms with duality gives $\nrd(a_p)=\overline{\nrd(a_p)}$. At $p\nmid r$, we use the identity map. Hence $u=(\nrd(a_p))_p$ satisfies $u=\bar u$, and \cref{lem:double-genus} gives $\hmod(L_r)\cong\hmod(L_1)$.
\end{proof}

\subsection{Proof of \texorpdfstring{\cref{thm:finite-all}}{Theorem B}}\label{sec:proof-B}

We first establish free summand cancellation (FSC) for two groups of order eight. The case $C_4\times C_2$ completes the group families needed for \cref{thm:finite-all}.

\begin{thm}\label{thm:small-free}
Let $G=C_4\times C_2$ or $C_2^3$. For any $\Z G$-lattices $P,Q$,
\[
 P\oplus\Z G\cong Q\oplus\Z G\quad\Longrightarrow\quad P\cong Q.
\]
Thus $\Z G$ has FSC, but it does not have TFC.
\end{thm}

\begin{proof}
Failure of TFC follows from \cite[Corollary~5.2]{Sw88}. Let $P,Q$ be $\Z G$-lattices with $P\oplus\Z G\cong Q\oplus\Z G$. The case $P=0$ follows by comparing ranks. Otherwise, local cancellation puts $P$ and $Q$ in the same genus. Since $G$ is abelian, $E=\End_{\Q G}(P_{\Q})$ is a product of matrix algebras over fields. At every odd prime $\ell$, $\Z_\ell G$ is a product of discrete valuation rings, so $\nu_\ell(P_\ell)$ contains all units of maximal orders on the support of $P$. If $P_2$ contains a free $R=\Z_2G$ summand, Lemma~\ref{lem:local-norm-sums} and Proposition~\ref{prop:norm-cancellation} give the result.

Otherwise let $N=\sum_{g\in G}g$. By Swan's decomposition lemma \cite[Lemma~5.6]{Sw88}, we can write $P$ as a direct sum of a projective $\Z G$-lattice and a lattice over $\Z G+\Z(N/2)$. The projective summand vanishes, since otherwise its completion at two would be free (since $\Z_2G$ is local). Thus $P_2$ is a lattice over $R'=R+\Z_2(N/2)$.

Each factor of $\Q_2G$ is $\Q_2$ or $\Q_2(i)$. Write $\mathcal O$ for the maximal order of $\Q_2G$ and $\mathfrak a$ for the ideal of $\mathcal O$ with component $4\Z_2$ on each rational factor and $2(1+i)\Z_2[i]$ on each factor $\Q_2(i)$. Any element of $\Q_2G$ whose group ring coefficients are all $\pm1/2$ differs from $N/2$ by an element of $R$. In particular $R'$ contains $4e_\chi$ for each rational character $\chi$. For $G=\langle x\rangle\times\langle y\rangle=C_4\times C_2$, the idempotents corresponding to factors $\Q_2(i)$ are $e_\pm=(1-x^2)(1\pm y)/4$. The coefficients of $2(1+x)e_\pm$ are again all $\pm1/2$, and $x$ acts as $i$ on these factors. These elements generate $\mathfrak a$ as an $R$-module, so $\mathfrak a\subseteq R'$.

Let $L=\mathcal O P_2$. Since $\mathfrak a L\subseteq P_2$, every unit of the maximal order on the support of $P$ that is congruent to $1$ modulo $\mathfrak a$ is realised as a reduced norm by scaling one basis vector in each supported factor of $L$. The automorphism and its inverse preserve $P_2$. On each $\Q_2$-factor, the global units $\Z^\times=\{\pm1\}$ represent both classes in $(\Z_2/4\Z_2)^\times$.
On each $\Q_2(i)$-factor, $2(1+i)\Z_2[i]=(1+i)^3\Z_2[i]$. The global units $\Z[i]^\times=\{\pm1,\pm i\}$ represent all four classes in $(\Z_2[i]/(1+i)^3\Z_2[i])^\times$, since the difference of any two distinct such units has $(1+i)$-adic valuation one or two.
Thus, for every unit $u$ of the maximal order on the support of $P$, we can choose a tuple $v$ of global units so that $uv^{-1}\equiv1\bmod\mathfrak a$ on the support. The argument above gives $uv^{-1}\in\nu_2(P_2)$. The fields $\Q$ and $\Q(i)$ have class number one, so every finite idele is a global element times units of maximal orders.
For any tuple of local units $(u_p)_p$, choose $v$ as above for $u_2$.
Then $u_pv^{-1}\in\nu_p(P_p)$ for every $p$, since at odd primes all
units of the maximal order are local norms. As $\nrd(E^\times)=K^\times$,
\eqref{eq:normclass} gives $\mathcal G(P)=1$, hence $P\cong Q$.
\end{proof}

It was shown by Hambleton--Kreck \cite[Proposition~2.4]{HK88} that, if $P = \pi_2(M)$ is the second homotopy module of a closed oriented $4$-manifold with finite fundamental group then $P$ is stably an extension of $L^*$ by $L$ where $L$ is a $\Z G$-lattice which fits into an exact sequence
\[
 0\longrightarrow L\longrightarrow F_2\longrightarrow F_1\longrightarrow F_0\longrightarrow\Z\longrightarrow0
\]
with the $F_i$ finitely generated free $\Z G$-modules. We refer to such a module $L$ as a \textit{third syzygy}.

Recall also that, if $P = \pi_2(M)$, then $P$ admits a hermitian form $h_P : P \times P \to \Z G$, given by the equivariant intersection form of the $4$-manifold. Furthermore, if $M$ and $N$ are stably diffeomorphic closed oriented smooth $4$-manifolds with $P = \pi_2(M)$ and $Q = \pi_2(N)$ then, after
identifying their fundamental groups and choosing compatible
orientations, the corresponding hermitian forms $h_P$ and $h_Q$ are related via an isometry
\[
 (P,h_P)\perp H(\Z G)^n \cong (Q,h_Q)\perp H(\Z G)^n
\]
where $H(\Z G)$ denotes the standard hyperbolic form on $\Z G^2$ (see, for example, \cite[\S 2]{HK93}).

We will now establish a cancellation theorem for $\Z G$-lattices satisfying these conditions.

\begin{thm}\label{thm:extension-cancellation}
Let $G$ be a finite group and let $P$ and $Q$ be $\Z G$-lattices. Suppose that, for some third syzygy $L$ of $\Z$ and some $a\ge0$, there is an exact sequence
\[
 0\longrightarrow L\longrightarrow P\oplus\Z G^a\longrightarrow L^*\longrightarrow0.
\]
Then $P\cong Q$ in either of the following cases.
\begin{clist}{(i)}
\item $G$ has periodic cohomology and $P\oplus\Z G^r\cong Q\oplus\Z G^r$ for some $r\ge0$.
\item $G$ has a cyclic normal subgroup of odd order with ambivalent $2$-group quotient, and $P$ and $Q$ carry hermitian forms which become isometric after adding the same number of standard hyperbolic forms.
\end{clist}
\end{thm}

\begin{proof}
\emph{Reduction to local norms.}
In either case $P$ and $Q$ become isomorphic after adding the same free summand. If $P=0$, comparing ranks gives $Q=0$. Otherwise local cancellation puts them in the same genus.

Fix a third syzygy $L$ and a partial free resolution as above. Write $d=\operatorname{rank}_{\Z G}F_2-\operatorname{rank}_{\Z G}F_1+\operatorname{rank}_{\Z G}F_0$. In the rational representation ring $[L_{\Q}]=d[\Q G]-[\Q]$. Since the regular and trivial representations are self-dual, the extension gives
\begin{equation}\label{eq:finite-rational-module}
 [P_{\Q}]=\chi[\Q G]-2[\Q],\qquad
 \chi=2d-a=2+\dim_{\Q}(P_{\Q})^G.
\end{equation}
In particular $\chi\ge2$. Each nontrivial simple $\Q G$-module occurs at least twice in $P_{\Q}$, and the endomorphism algebra on the trivial component, if present, is split. Thus $E=\End_{\Q G}(P_{\Q})$ satisfies the Eichler condition.

By Proposition~\ref{prop:norm-cancellation} and Lemma~\ref{lem:local-norm-sums}, free summands cancel from $P$ if
\begin{equation}\label{eq:local-norm-cancellation}
 \pi_{P,p}\bigl(\nu_p(\Z_pG)\bigr)\subseteq\nu_p(P_p)
\end{equation}
at every prime $p$.

\emph{Cyclic and generalised quaternion Sylow subgroups.}
We prove \eqref{eq:local-norm-cancellation} whenever a Sylow $p$-subgroup of $G$ is cyclic or generalised quaternion. For $p\nmid|G|$, the order $\Z_pG$ is maximal \cite[Proposition~27.1]{CR81}.
By \cite[Theorems~26.20(i) and~26.24(iii)]{CR81}, each simple component of $\Q_pG$ supports a unique indecomposable $\Z_pG$-lattice up to isomorphism. Thus every indecomposable projective on the rational support of $P_p$ occurs as a direct summand of $P_p$ itself, and Lemma~\ref{lem:local-norm-sums} applies.

Suppose $p\mid|G|$, and let $D$ be a Sylow $p$-subgroup. Since $D$ is cyclic or generalised quaternion, its period divides four \cite[XII, proof of Theorem~11.6]{CE56}. Thus Tate periodicity gives $H_4(D;\Z)=\HH^{-1}(D;\Z)=0$ \cite[VI, \S9]{Br82}.
On $H_4(G;\Z)$, transfer to $D$ followed by the map induced by inclusion is multiplication by $[G:D]$ \cite[III, Proposition~9.5]{Br82}.
Since this composite factors through $H_4(D;\Z)=0$, the index $[G:D]$ annihilates the extension class under the isomorphism $\Ext^1_{\Z G}(L^*,L)\cong H_4(G;\Z)$ \cite[Proposition~2.4(iii)]{HK88}.
As $[G:D]$ is a unit in $\Z_p$, the completed extension class is zero, so the extension splits over $\Z_pG$.
Let $R=\Z_pG$, $k=\F_p$, and take a minimal projective resolution
\[
 0\longrightarrow J\longrightarrow E_2\longrightarrow E_1\longrightarrow E_0\longrightarrow\Z_p\longrightarrow0.
\]
Since $L_p$ and $J$ are both third syzygies of $\Z_p$ over $R$, Schanuel's lemma identifies them up to projective summands. The same holds for their duals, so the splitting $P_p\oplus R^a\cong L_p\oplus L_p^*$ identifies $P_p$ with $J\oplus J^*$ up to projective summands.
Reduction modulo $p$ gives a minimal resolution over the self-injective algebra $kG$, whose syzygies have no projective summands. Thus $J,J^*$ have none, and Krull--Schmidt gives $P_p\cong J\oplus J^*\oplus U$ with $U$ projective.

Suppose the $p$-period of $G$ divides four and let $I=\ker(\eps:\Z G\to\Z)$.
The dual augmentation sequence
\[
  0\longrightarrow\Z_p\xrightarrow{\,1\mapsto\sum_{g\in G}g\,}R
  \longrightarrow I_p^*\longrightarrow0
\]
identifies $I_p^*$ with $\Omega_R^{-1}\Z_p$ and hence, by four-periodicity, with $\Omega_R^3\Z_p=J$ up to projective summands.
By Krull--Schmidt, $P_p$ and $I_p\oplus I_p^*\oplus R^{\chi-2}$ have the same nonprojective indecomposable summands with the same multiplicities. Their direct sum is $J\oplus J^*$. Subtracting its character from the common rational character in \eqref{eq:finite-rational-module}, we see that the remaining projective parts have the same rational character. Swan's theorem \cite[Corollary~6.4]{Sw60} therefore yields
\[
 P_p\cong I_p\oplus I_p^*\oplus R^{\chi-2}.
\]
For $\chi>2$ there is a free summand. For $\chi=2$, right multiplication on $I_p$ by each unit of $R$ realises its regular reduced norm projected away from the trivial component. Both give \eqref{eq:local-norm-cancellation}.

\emph{The case of finite $p$-period greater than four.}
Here $p$ is odd and $D$ is cyclic \cite[Theorem~1]{Sw60p}. By Swan's formula \cite[Theorem~2]{Sw60p}, the $p$-period of $G$ is $2e$, where $e=|N_G(D)/C_G(D)|$, so $e>2$.
In case~\textup{(ii)}, $D$ would be the unique Sylow $p$-subgroup of the cyclic normal subgroup $C$, so $D$ would be normal in $G$ and conjugation on $D$ would factor through $G/C$.
Since $G/C$ is ambivalent, every element of its image in the cyclic group $\Aut(D)$ would equal its inverse. This would give $e\le 2$, a contradiction.
Thus case~\textup{(ii)} cannot occur here. We now assume hypothesis~\textup{(i)}.

For a simple $kG$-module $S$, let $D_S=\End_{kG}(S)$, $d(S)=\dim_{D_S}S$ and $h(S)=\dim_{D_S}H^2(G;S)$. Restriction gives a $D_S$-linear injection $H^2(G;S)\to H^2(D;S)$ \cite[III, Proposition~9.5]{Br82}. Since $D$ is cyclic, $H^2(D;S)\cong S^D/N_DS$, where $N_D=\sum_{g\in D}g$. Hence
\[
 h(S)\le \dim_{D_S}H^2(D;S)=\dim_{D_S}(S^D/N_DS)
 \le \dim_{D_S}S^D\le \dim_{D_S}S=d(S).
\]
If $h(S)=d(S)$, equality in the first inequality makes restriction surjective, and equality in the last gives $S^D=S$. Thus $D$ acts trivially on $S$, so the target of restriction is naturally $H^2(D;k)\otimes_k S$.
Let $\lambda$ be the character of $N_G(D)$ on $H^2(D;k)$. The image of restriction is fixed by $N_G(D)$, so in the equality case the normaliser acts on $S$ by $\lambda^{-1}$.
Write $|D|=p^m$. By \cite[Lemma~3]{Sw60p}, the pullback of $x\mapsto x^r$ acts by $r$ on $H^2(D;\Z)$, hence by $r\bmod p$ on $H^2(D;k)\cong H^2(D;\Z)\otimes_{\Z}k$.
Reduction $(\Z/p^m\Z)^\times\to k^\times$ has $p$-group kernel, so it is injective on the conjugation image $N_G(D)/C_G(D)$, whose order is prime to $p$ since $D\le C_G(D)$. Thus $\lambda$ has order $e>2$, and $h(k)=0$.
If equality holds for $S$, the normaliser acts on the contragredient dual $S^*$ by $\lambda$. Equality for $S^*$ would require it to act by $\lambda^{-1}$ as well, forcing $\lambda^2=1$, contrary to $e>2$. Hence  $h(S)+h(S^*)<2d(S)$. 
 
The resolution and rational character of $P_p$ give, by Swan's injectivity theorem,
\[
 [U]=\chi[R]-2[E_0]+[E_1]+[E_1^*]-[E_2]-[E_2^*].
\]
The projective cover of $S$ occurs $d(S)$ times in $R$ and $h(S)$ times in $E_2$, since applying $\Hom_{kG}(-,S)$ to the minimal resolution modulo $p$ gives zero differentials. In $2[R]-2[E_0]-[E_2]-[E_2^*]$, the multiplicity of this cover is therefore zero for $S=k$ and $2d(S)-h(S)-h(S^*)>0$ otherwise, since $E_0$ is the self-dual cover of $k$. Thus $U$ contains $R^{\chi-2}\oplus E_1$, and every indecomposable projective other than $E_0$ occurs in $U$. If $\chi>2$, this gives the required free summand of $P_p$.

Suppose $\chi=2$. By Lemma~\ref{lem:local-norm-sums}, $\nu_p(R)$ is the product of the norm groups of its distinct indecomposable projective summands. Every indecomposable projective other than $E_0$ already occurs as a direct summand of $P_p$. Thus, to prove \eqref{eq:local-norm-cancellation}, it suffices to show that $\pi_{P,p}(\nu_p(E_0))\subseteq\nu_p(P_p)$.
Write $\eps:E_0\to\Z_p$ for the projective cover. Since $\Hom_R(E_0,\Z_p)=\Z_p\eps$, every $u_0\in\Aut_R(E_0)$ satisfies $\eps u_0=c\eps$ for some $c\in\Z_p^\times$ and hence preserves $\ker\eps$.
Lift the restriction of $u_0$ to $\ker\eps$ through its projective cover $E_1$. The lift $u_1$ induces an isomorphism on the top, so it is an automorphism by Nakayama's lemma. Repeating this argument at the next kernel gives an automorphism $u_2$ of $E_2$, which restricts to an automorphism $u_3$ of $J$.
Since $\chi=2$, the projection $\pi_{P,p}$ omits the trivial component by \eqref{eq:finite-rational-module}. Hence multiplicativity along the rational resolution gives
\[
 \pi_{P,p}(\nrd(u_0))=\nrd(u_3)\nrd(u_1)\nrd(u_2)^{-1}\in\nu_p(P_p),
\]
since $J$ and $E_1$ are summands of $P_p$ and every indecomposable summand of $E_2$ occurs in $U$ because $h(k)=0$. This proves \eqref{eq:local-norm-cancellation}.

In case \textup{(i)}, every Sylow subgroup is cyclic or generalised quaternion \cite{Sw60p}, so Proposition~\ref{prop:norm-cancellation} completes the proof. In case \textup{(ii)}, every odd Sylow subgroup lies in the given cyclic normal subgroup, so only the prime two remains.

\emph{The prime two in case \textup{(ii)}.}
Let $C$ be the given cyclic normal subgroup and $T=G/C$. Put
\[
 e=|C|^{-1}\sum_{c\in C}c,\qquad
 R_0=e\Z_2G\cong\Z_2T,\qquad R_1=(1-e)\Z_2G .
\]
Since $e$ is a central idempotent, $R=\Z_2G\cong R_0\times R_1$. Reduced norms are computed componentwise, so to prove \eqref{eq:local-norm-cancellation} at two it suffices to verify the inclusion over $R_0$ and $R_1$ separately.
Let $J_1=(1-e)L_2$ and $P_1=(1-e)P_2$. Applying $1-e$ to the completed resolution of $\Z$ kills $\Z_2$. Splitting successively from the right shows that $J_1$ is projective. Its dual $J_1^*\cong(1-e)(L^*)_2$ is also projective since $\bar e=e$. The completed extension therefore splits on $R_1$, giving
\[
 P_1\oplus R_1^a\cong J_1\oplus J_1^*.
\]
Thus, $P_1$ is projective. Its rational character is that of $R_1^\chi$, so $P_1\cong R_1^\chi$ by Swan's theorem \cite{Sw60}. This gives \eqref{eq:local-norm-cancellation} on $R_1$.

Write $P_0=eP_2$ and $Q_0=eQ_2$. If $P_0$ contains a free $R_0$-summand, \eqref{eq:local-norm-cancellation} holds on both factors at two and the proof is complete. Suppose then that $P_0$ has no free summand. By local cancellation $Q_0\cong P_0$, so every $R_0$-linear map $Q_0\to R_0$ has image in the radical of the local ring $R_0$.

Choose a rational isometry between the forms on $P_{\Q}$ and $Q_{\Q}$ by Witt cancellation \cite[Corollary~5.3]{BFM14}. Let $n$ be the number of added hyperbolic forms. If $n=0$, the given isometry already gives $P\cong Q$, so assume $n\geq 1$. In these coordinates let $f$ be the stable isometry and write $\eta=\pi_P(\nrd(f))\in\nrd(E^\times)$ as in Proposition~\ref{prop:norm-cancellation}.

We write $H=H(R_0^n)$ and $f_0:P_0\perp H\to Q_0\perp H$ for the $R_0$-component of $f$. We construct an isometry $s$ of $Q_0\perp H$ such that $s\circ f_0=a_0\perp\id_H$ for an isometry $a_0:P_0\to Q_0$. To control the reduced norm of $a_0$, we track the reduced norm of each elementary step as we send the image of each hyperbolic pair to the corresponding standard pair.

We write $h$ for the form on $Q_0\perp H$. If $(x,y)$ is the image under $f_0$ of the first standard pair, the $R_0$-linear map $h(-,y)$ sends $Q_0$ into the radical, so $h(x,y)=1$ forces one hyperbolic coordinate of $x$ to be a unit. We permute the planes, flip a pair if necessary, and scale its two generators by $u$ and $\bar u^{-1}$ for a unit $u$. Then $x=v+z+bw$, with $(v,w)$ the first standard pair and $z$ in its orthogonal complement. Since $x$ is isotropic, $b+\bar b=-h(z,z)$.

With $h$ linear in its first variable, the triangular isometry
\[
 v\longmapsto v+z+bw,\qquad w\longmapsto w,\qquad
 t\longmapsto t-h(t,z)w\quad(t\perp v,w)
\]
has reduced norm one. We apply its inverse to send $x$ to $v$. The image of $y$ then has the form $w+z'+b'v$, since it pairs to $1$ with $v$. Interchanging $v$ and $w$ in the construction, we send this vector to $w$ while fixing $v$. We can now split off this hyperbolic plane and repeat the argument on its orthogonal complement. We take $s$ to be the composite of these isometries.

Permuting whole planes has reduced norm one, and scaling a pair has norm $\nrd(u)/\overline{\nrd(u)}=1$ since $T$ is ambivalent. A flip has the norm of $-\id$ on one regular summand, its rational $(-1)$-eigenspace.
Set $\delta=\pi_P(\nrd(-\id_{\Q G}))$. By \eqref{eq:finite-rational-module}, the regular $\Q G$-module restricted to the rational support of $P$ occurs as a direct summand of $P_{\Q}$. Acting by $-1$ on this summand and by the identity on a complement gives an automorphism with reduced norm $\delta$. Thus, $\delta\in\nrd(E^\times)$.
Since $s\circ f_0=a_0\perp\id_H$, we obtain $\nrd(a_0)=(\eta\delta^j)_0$ for some $j\in\{0,1\}$, where the subscript denotes the $R_0$-component.

Choose an $R_1$-isomorphism $a_1:P_1\to(1-e)Q_2$. Comparing $a_1\oplus\id$ with the $R_1$-component of $f$ and using \eqref{eq:local-norm-cancellation} gives $\nrd(a_1)/\eta_1\in\nu_2(P_1)$. Since $P_1\cong R_1^\chi$, the norm $\delta_1$ of $-\id$ on one regular summand also lies in $\nu_2(P_1)$. Together with $\nrd(a_0)=(\eta\delta^j)_0$, this gives
\[
 a_2=a_0\oplus a_1:P_2\xrightarrow{\ \cong\ }Q_2,\qquad
 \frac{\nrd(a_2)}{\eta\delta^j}\in\nu_2(P_2).
\]
At each odd prime $\ell$, choose an isomorphism $a_\ell:P_\ell\to Q_\ell$. Comparing $a_\ell\oplus\id$ with $f_\ell$ gives $\nrd(a_\ell)/\eta\in\nu_\ell(P_\ell)$ by \eqref{eq:local-norm-cancellation}. Applying the same inclusion to $-\id$ on $\Z_\ell G$ gives $\delta\in\nu_\ell(P_\ell)$, so $\nrd(a_\ell)/(\eta\delta^j)\in\nu_\ell(P_\ell)$ as well. Dividing the norm idele by the global norm $\eta\delta^j$ therefore gives a local automorphism norm at every prime. Formula~\eqref{eq:normclass} gives $P\cong Q$, proving \textup{(ii)}.
\end{proof}

\begin{proof}[Proof of \cref{thm:finite-all}]
Use the stable diffeomorphism to identify the fundamental groups and to choose compatible orientations. For $P=\pi_2(M)$ the extension in Theorem~\ref{thm:extension-cancellation} is precisely Hambleton--Kreck's Proposition~2.4(i) \cite{HK88}. Connected sum with $S^2\times S^2$ adds a standard hyperbolic plane to the equivariant intersection form \cite[\S2]{HK93}. Thus the forms on $P$ and $Q=\pi_2(N)$ are stably isometric, and forgetting the forms gives a stable isomorphism of the modules. 

Theorem~\ref{thm:extension-cancellation} proves \textup{(i)} and \textup{(ii)}. 
For \textup{(iii)}, the cyclic groups, $S_3$, and $Q_8$ are covered by \textup{(i)}, while $C_2^2$, $C_2^3$, and the dihedral group $D_8$ of order eight are covered by \textup{(ii)}. The only remaining group of order at most eight is $C_4\times C_2$, which is covered by Theorem~\ref{thm:small-free}.
Forgetting the fixed identifications of the fundamental groups gives the stated $\Aut(G)$-isomorphism.
\end{proof}

\begin{corollary}\label{cor:periodic-euler}
Let $X$ and $Y$ be finite oriented Poincar\'e $4$-complexes whose fundamental groups are identified with the same finite group $G$ of periodic cohomology. If $\chi(X)=\chi(Y)$, then $\pi_2(X)\cong\pi_2(Y)$ as $\Z G$-modules.
\end{corollary}

\begin{proof}
The transfer argument in the proof of Theorem~\ref{thm:extension-cancellation} gives $H_4(G;\Z)=0$, so the Hambleton--Kreck extensions for $X$ and $Y$ split integrally. The two third syzygies become isomorphic after adding free summands by Schanuel's lemma. Taking duals and using the splittings gives $\pi_2(X)\oplus \Z G^m\cong\pi_2(Y)\oplus\Z G^n$ for some $m,n$. Since $\chi(X)=\chi(Y)$, the universal covers have the same second Betti number and hence $m=n$. Now apply Theorem~\ref{thm:extension-cancellation}\textup{(i)}.
\end{proof}

\subsection{Proof of \texorpdfstring{\cref{thm:finite-doubles}}{Theorem C}}

\leavevmode Fix identifications $\pi_1(X)\cong G$ and $\pi_1(Y)\cong G$,
and let $L=\pi_2(X)$ and $L'=\pi_2(Y)$.
With the induced markings on the doubles,
\cite[Proposition 5.2]{HN26} gives
\begin{equation}\label{eq:finite-double}
 \pi_2(M(X))\cong\hmod(L),\qquad \pi_2(M(Y))\cong\hmod(L').
\end{equation}

In case (i), since $\chi(X)=\chi(Y)$, generalised
Schanuel \cite[VIII, Lemma 4.4]{Br82}, applied to the
augmented cellular sequences, gives
$L\oplus\Z G^a\cong L'\oplus\Z G^a$ for some $a\ge0$.
Local cancellation \cite[Theorem 20.6.2]{Vo21} then gives
$L_p\cong L'_p$ for every prime, so $L$ and $L'$ lie
in the same genus.

If $L=0$, then $L'=0$ and the conclusion is immediate.
We may therefore assume $L\ne0$.

The centre of $\Q G$ is spanned by conjugacy-class sums.
Ambivalence implies that inversion fixes each such sum,
so the involution on the supported centre in
\cref{lem:double-genus} is trivial, also after completion.
Thus $u\bar u^{-1}=1$ in the notation of that lemma,
which gives $\hmod(L)\cong\hmod(L')$ with the chosen
markings fixed.

In case (ii), apply \cref{lem:power} at each commutator
twist, using the inverse isomorphism when a step is traversed
backwards. Composing these isomorphisms identifies the doubled
modules of the endpoint presentation complexes, preserving
the common marking by the displayed generators.
The endpoint homotopy equivalences give $\Aut(G)$-isomorphisms
with $\hmod(L)$ and $\hmod(L')$, since duality commutes with
changes of marking. In either case, \eqref{eq:finite-double}
proves the result.

\subsection{Proof of \texorpdfstring{\cref{cor:abelian}}{Corollary D}}

The trivial group is immediate. We use Browning's normal form \cite[Theorem 1.7]{Bro79}, together with the classification stated in \cite[Theorem 8.1]{HN26}. Above the minimum Euler characteristic there is one homotopy type at each Euler characteristic. At the minimum, for
$G=C_{m_1}\times\cdots\times C_{m_d}$ with $1<m_1\mid\cdots\mid m_d$ and $d\ge2$, every complex is homotopy equivalent to the presentation complex
\[
 \left\langle x_1,\ldots,x_d\ \middle|\
 x_i^{m_i},\ [x_1^r,x_2],\ [x_i,x_j]\ ((i,j)\ne(1,2),\ i<j)
 \right\rangle,
\]
where $(r,m_1)=1$. These are the presentation complexes of Sieradski--Dyer \cite[Proposition 6]{SiDy79}. Each differs from the $r=1$ presentation by a commutator twist, so \cref{thm:finite-doubles}(ii) identifies their doubled modules. The cyclic case has a unique homotopy type at each Euler characteristic \cite[Theorem A]{DySi73full}. Since the normal forms are obtained up to unmarked homotopy equivalence, the conclusion is an $\Aut(G)$-isomorphism.

\begin{remark}\label{rem:ambivalent}
Springer \cite[Theorem 8.7]{Sp74}, extending Carter's result for Weyl groups \cite{Ca72}, proves that every element of a finite Coxeter group is a product of two elements whose squares are the identity. If $g=ab$ with $a^2=b^2=1$, then $aga^{-1}=g^{-1}$, so these groups are ambivalent. This includes all dihedral and symmetric groups. Also, $\operatorname{Dih}(A)=A\rtimes_{-1}C_2$ is ambivalent for every finite abelian group $A$: the generator of $C_2$ inverts $A$, and every element outside $A$ is an involution. Direct products preserve ambivalence. Each of these families therefore satisfies \cref{thm:finite-doubles}(i).
\end{remark}

\begin{remark}\label{rem:quadratic-bias}
Within each $4$-dimensional family distinguished by the quadratic bias invariant in \cite[Section 9]{HN26}, the second homotopy modules are $\Aut(G)$-isomorphic as $\Z G$-modules, where $G$ is the common fundamental group. For the abelian families this follows from \cref{cor:abelian}. The non-abelian family has group $Q_8\times(\Z/p)^3$, where $p$ is prime and $p\equiv1\pmod8$. Its presentations contain $C^p$ and differ only by replacing $[A,C]$ with $[A,C^r]$, where $(r,p)=1$ \cite[Lemma 9.6]{HN26}. Inverting these relators gives $[C,A]\mapsto[C^r,A]$, so \cref{thm:finite-doubles}(ii) applies. Thus the quadratic bias invariant distinguishes homotopy types even when their second homotopy modules agree.
\end{remark}

\section{Homotopy equivalent doubles}\label{sec:noninjective}

We prove \cref{thm:noninjective}. The argument compares the full quadratic $2$-types of the two doubles; the cancellation theorem alone does not give their homotopy equivalence.

\subsection{Norm extensions}

For the complexes $X,Y$ of \cref{thm:noninjective}, let $G=Q_{28}$, $N=\sum_{g\in G}g$, $L_0=\pi_2(X)$ and $L_1=\pi_2(Y)$. We use the duality conventions of \cref{sec:finite}. The hyperbolic form $H(L)$ is the $G$-invariant integral symmetric form on $L\oplus L^*$ given by
\[
 ((x,f),(y,h))\longmapsto f(y)+h(x).
\]
Such forms are equivalent to $\Z G$-valued hermitian forms for the involution $\overline g=g^{-1}$: the coefficient of $1$ gives the integral form, and the other coefficients are recovered by translating an argument. We use this correspondence throughout the section. The underlying module of $H(L)$ is $\hmod(L)$.

Fix the markings of the fundamental groups supplied by the generators $x,y$. Let $\kappa_i\in H^3(G;L_i)$ be the first $k$-invariant. The standard complex has
\begin{equation}\label{c:eq:standard-extension}
 L_0\cong \Z G/\Z N,\qquad
 0\longrightarrow\Z\xrightarrow{1\mapsto N}\Z G\longrightarrow L_0\longrightarrow0;
\end{equation}
see \cite[Section 2]{MP21}. Its connecting homomorphism sends $\kappa_0$ to a generator $u\in H^4(G;\Z)\cong\Z/28$. Indeed, dimension shifting along the augmented cellular complex identifies $\kappa_0$ with $1\in\widehat H^0(G;\Z)=\Z/28$, and the connecting map in \eqref{c:eq:standard-extension} is an isomorphism in positive degrees. This fixes the choice of $u$ and the sign convention for connecting maps.

\begin{lemma}\label{c:lem:extension}
There is a stably free rank-one $\Z G$-module $P$, with $P\oplus \Z G\cong \Z G^2$, and an exact sequence
\begin{equation}\label{c:eq:P-extension}
 0\longrightarrow\Z\xrightarrow{1\mapsto e}P\longrightarrow L_1\longrightarrow0
\end{equation}
whose connecting map sends $\kappa_1$ to $u$. Moreover, $e$ generates $P^G=NP$ and is primitive as an element of the underlying abelian group.
\end{lemma}
\begin{proof}
The simple homotopy equivalence $Y\vee S^2\simeq_s X\vee S^2$ can be chosen over $G$. Indeed, the presentations have one common relator, and adjoining the other second relator to either presentation gives the same three-relator presentation. The added relator is null-homotopic, so its attachment is simple homotopy equivalent to wedging on $S^2$. This construction induces the identity on the labelled generators. Hence there is an isomorphism
\[
 t:L_1\oplus \Z G\xrightarrow{\cong}L_0\oplus \Z G,
 \qquad t_*(\kappa_1,0)=(\kappa_0,0).
\]
Pull back along $t$ the exact sequence
\[
 0\longrightarrow\Z\xrightarrow{1\mapsto(N,0)}\Z G^2
 \longrightarrow L_0\oplus \Z G\longrightarrow0.
\]
The middle module is isomorphic to $\Z G^2$. Its surjection onto the free summand of $L_1\oplus \Z G$ splits; call its kernel $P$. This gives \eqref{c:eq:P-extension}, and naturality gives the assertion about $u$.

Since $L_1$ is a lattice, $e$ is primitive. Projectivity implies that the norm map $P_G\to P^G$ is an isomorphism. Both groups are infinite cyclic, as follows also from $P\oplus \Z G\cong \Z G^2$. The non-zero invariant primitive element $e$ therefore generates $P^G=NP$.
\end{proof}

\subsection{Hyperbolic forms and the \texorpdfstring{$k$}{k}-invariant}

For a closed oriented $4$-manifold $M$, its quadratic $2$-type consists of $\pi_1(M)$, the $\pi_1(M)$-module $\pi_2(M)$, its first $k$-invariant, and its equivariant intersection form. An isometry must preserve all this data, including the $k$-invariant.

Let $X$ be a finite $2$-complex with finite fundamental group, and let $L=\pi_2(X)$ and $\kappa=k_X$. The thickening model for its standard double has a retraction to $X$ and a homotopy section. Poincar\'e--Lefschetz duality identifies the kernel on $\pi_2$ with $L^*$. With the splitting given by the section, the quadratic $2$-type is
\[ 
\bigl(G,L\oplus L^*,(\kappa,0),\Met(L,\phi)\bigr),\qquad
 \Met(L,\phi)=\begin{pmatrix}\phi&1\\1&0\end{pmatrix},
\]
where $\phi:L\to L^*$ is symmetric. This is the metabolic description of \cite[Propositions II.2 and II.4]{KS84}. Naturality with respect to the section gives the $k$-invariant $(\kappa,0)$.

\begin{lemma}\label{c:lem:hyperbolic}
If the finite group $G$ has a cyclic Sylow $2$-subgroup, then the quadratic $2$-type of a standard double $M(X)$ is isometric to
\[
 \bigl(G,L\oplus L^*,(\kappa,0),H(L)\bigr).
\]
The isometry can be chosen to induce the identity on the quotient $L$ of the retraction.
\end{lemma}
\begin{proof}
It suffices to find a $G$-homomorphism $b:L\to L^*$ such that
\begin{equation}\label{c:eq:good-b}
 b+b^*=\phi,\qquad b_*\kappa=0.
\end{equation}
The shear $(x,f)\mapsto(x,f+bx)$ then gives the required isometry.

\emph{The cyclic case.}
By Dyer--Sieradski \cite[Theorem and Corollary]{DySi73}, a finite $2$-complex with cyclic fundamental group is simple homotopy equivalent to the standard one-relator complex wedged with $2$-spheres. The one-relator complex is a spine of a punctured lens space $V$. Suciu's calculation of its spin gives the hyperbolic form and $k$-invariant $(\kappa,0)$ \cite[Section 4]{Suc88}; the sphere summands add hyperbolic free summands.

Let $W_0$ be the boundary connected sum of $V\times D^2$ with the required copies of $S^2\times D^3$, and let $W$ be the given thickening. Both are parallelisable. Phillips's theorem \cite[Theorem A]{Ph67} and general position on a $2$-dimensional spine realise the chosen homotopy equivalence by an orientation-preserving embedding $W_0\hookrightarrow\operatorname{int}(W)$: first immerse the interiors, then restrict to a small regular neighbourhood of the embedded spine. The complement is an $h$-cobordism, by van Kampen's theorem, excision and duality; compare with \cite[p.~15]{KS84}.

The resulting boundary equivalence commutes up to homotopy with the maps to the spines. After identifying quotient modules, it induces the identity on $L$; the evaluation pairing then forces the identity on $L^*$. Its isometry therefore preserves the $k$-invariant, and gives $b$ satisfying \eqref{c:eq:good-b}.

\emph{Transfer from a Sylow subgroup.}
Let $C\le G$ be a cyclic Sylow $2$-subgroup and let $m=[G:C]$. The corresponding cover of the thickening is stably parallelisable. The cyclic case gives a $\Z C$-homomorphism $b_C$ satisfying \eqref{c:eq:good-b} after restriction. Its transfer
\[
 b_m=\sum_{g\in G/C}g b_C g^{-1}
\]
satisfies $b_m+b_m^*=m\phi$ and
$(b_m)_*\kappa = \cor_C^G\bigl((b_C)_*\res_C^G\kappa\bigr)=0$.
The latter equality is the transfer projection formula for the evaluation pairing $\Hom_{\Z}(L,L^*)\otimes_{\Z}L\to L^*$.

\emph{Removing the odd factor.}
Reflection in the double commutes with the retraction, reverses orientation and preserves the $k$-invariant. Its action is therefore
\[
 J=\begin{pmatrix}1&0\\ c&-1\end{pmatrix},
 \qquad c+c^*=-2\phi,\qquad c_*\kappa=0.
\]
Since $m$ is odd, $b=b_m+\frac{m-1}{2}c$ satisfies \eqref{c:eq:good-b}.
\end{proof}

For $Q_{28}$ one can take $C=\langle y\rangle\cong C_4$, so $m=7$. In particular, the lemma applies to $X$ and $Y$.

\subsection{Integral isometries}

Here $G=Q_{28}$. For a rational prime $\ell$, we use the group rings $\Z_\ell G$ and $\Q_\ell G$, and write $P_\ell=\Z_\ell\otimes_{\Z}P$ for a lattice $P$. We write $H(P)$ for the standard hyperbolic form on $P \oplus P^*$.

\begin{prop}\label{c:prop:arithmetic}
Let $P\oplus \Z G\cong \Z G^2$, and choose a generator $e$ of $P^G$. There is an isometry
\[
 U:H(\Z G)\xrightarrow{\cong}H(P),\qquad U(N,0)=(e,0).
\]
Moreover, for any positive integer $a$, $U$ can be chosen so that, for each $\ell\mid a$, there is an isomorphism $f_\ell:\Z_\ell G\to P_\ell$, with $f_\ell(N)=e$, such that
\[
 \bigl(U_\ell-H(f_\ell)\bigr)H(\Z_\ell G)\subseteq aH(P_\ell).
\]
Here $H(f)=f\oplus(f^*)^{-1}$.
\end{prop}

\begin{proof}
\emph{Normalising the local comparisons.}
Choose an identification $P\otimes\Q\cong \Q G$ carrying $e$ to $N$, and regard $P$ as a lattice in $\Q G$. The norm identifies $P_G$ with $P^G$, so $\eps(P)=\Z$, where $\eps:\Q G\to\Q$ is augmentation. Locally $P_\ell=\Z_\ell G a_\ell$ for some $a_\ell\in (\Q_\ell G)^\times$. This follows from cancellation over the semilocal ring $\Z_\ell G$. Since $\eps(P_\ell)=\Z_\ell$, the generator can be chosen with $\eps(a_\ell)=1$.

Use row vectors for left modules, so module maps are right multiplication by matrices. A stable trivialisation $P\oplus \Z G\to \Z G^2$ extends to a matrix $B\in\GL_2(\Q G)$. By a change of basis in $\GL_2(\Z)$ on the target, arrange that $\eps(B)=I_2$. Then
\[
 C_\ell=\operatorname{diag}(a_\ell,1)B\in\GL_2(\Z_\ell G).
\]
Write $\nrd$ for the componentwise reduced norm, also on matrix algebras. Semilocal elimination \cite[proof of Lemma III.1.4]{We13}, applied to $\Q G$ and $\Z_\ell G$, gives $b\in (\Q G)^\times$ and $c_\ell\in (\Z_\ell G)^\times$ such that
\[
 \nrd(b)=\nrd(B),\qquad \nrd(c_\ell)=\nrd(C_\ell).
\]
Their augmentations are $1$, since the reduced norm on the trivial component is the ordinary determinant and $\eps(B)=\eps(C_\ell)=I_2$.

Replace $P$ by the isomorphic lattice $P'=Pb$. Right multiplication by $b$ preserves $N$. The local generators $a'_\ell=c_\ell^{-1}a_\ell b$ satisfy $P'_\ell=\Z_\ell G a'_\ell$ and $\nrd(a'_\ell)=1$.
Here left multiplication by $c_\ell^{-1}\in (\Z_\ell G)^\times$ preserves the left ideal $\Z_\ell G a_\ell b$. For all but finitely many $\ell$, take $a'_\ell=1$.

\emph{The rational factors.}
The rational algebra with involution has the decomposition
\[
 \Q G\cong\Q\times\Q\times\Q(i)\times M_2(F)\times D,
 \qquad F=\Q(\zeta_7+\zeta_7^{-1})
\]
(see, for example, the $n=7$ case of \cite[Section 6, p.~11848]{Ni24}). Here
$D=\Q[\zeta_{14},j]$ is a totally definite quaternion algebra. The involution $g\mapsto g^{-1}$ restricts to conjugation on $D$, since $x\mapsto\zeta_{14}$ and $y\mapsto j$.
We can identify the $M_2(F)$-factor so that the involution is
$A\mapsto B^{-1}A^tB$, where
$B=\left(\begin{smallmatrix}2&s\\s&2\end{smallmatrix}\right)$ and
$s=\zeta_7+\zeta_7^{-1}$.

\emph{Approximation.}
We use strong approximation on the two noncommutative factors to
find a single rational isometry carrying $H(\Z G)$ onto $H(P')$.
To this end, we will approximate the local isometries induced by multiplication with $a'_\ell$ at finitely many primes using strong approxmiation.
Let $f'_\ell\colon\Z_\ell G\to P'_\ell$ be right multiplication by $a'_\ell$. Then $H(f'_\ell)$ carries $H(\Z_\ell G)$ onto $H(P'_\ell)$ and is the identity on the three commutative factors, since $\nrd(a'_\ell)=1$.
On $M_2(F)$, viewing a $2\times2$ matrix $T$ over $M_2(F)$ as a
$4\times4$ matrix over $F$, the isometry condition becomes
$T^tQT=Q$, where
\[
 Q=\begin{pmatrix}0&B\\ B&0\end{pmatrix}.
\]
Both coordinate planes are totally isotropic, so the quadratic form
with Gram matrix $Q$ is split.
Since $\nrd(a'_\ell)=1$, its $M_2(F_\ell)$ component lies in
$\SL_2(F_\ell)$. This group is generated by elementary matrices,
each of which is a square. The induced isometry of $Q$ therefore
has trivial spinor norm and lies in the spinor kernel.
On $D$, the form is a hyperbolic quaternionic hermitian plane.
Both forms are isotropic at the real places, so strong
approximation applies to $\operatorname{Spin}(Q)$ on the $M_2(F)$-factor
and to the unitary group of the hyperbolic hermitian plane on the
$D$-factor (see \cite{Kn66}, \cite[Theorem 5.1.3]{Kir16}).
The norm-one group of $D$ is compact at every real place, which is why
we apply strong approximation after hyperbolising.

Choose a finite set $S$ of primes containing those dividing $a$ and
large enough that $H(\Z_\ell G)=H(P'_\ell)$ for every $\ell\notin S$.
At each prime outside $S$, require the assembled isometry to preserve
this common lattice. At primes in $S$,
lift the $M_2(F_\ell)$-component of $H(f'_\ell)$ to $\operatorname{Spin}(Q)(F_\ell)$.
We project the spin component back to $\operatorname{SO}(Q)$ and assemble
it with the quaternionic component and the identity on the commutative
factors. We impose the conditions on this assembled map: it must carry
the whole lattice $H(\Z_\ell G)$ onto $H(P'_\ell)$ and agree with
$H(f'_\ell)$ modulo $aH(P'_\ell)$ when $\ell\mid a$. Strong approximation gives
a rational isometry satisfying these open conditions.

The resulting isometry identifies $H(\Z G)$ with $H(P')$ integrally. It fixes $(N,0)$ because it acts identically on the trivial rational component. Transport it back along $P\to P'$, and use the transported local multiplication maps as the $f_\ell$.
\end{proof}

\subsection{Orthogonal reduction}

For a primitive invariant element $e\in P$, orthogonal reduction of $H(P)$ gives a canonical identification
\[
 (e,0)^\perp/\Z(e,0)
 = (P/\Z e)\oplus(P/\Z e)^*
 \cong H(P/\Z e)
\]
as forms. Indeed, $(e,0)^\perp=P\oplus\{f\in P^*:f(e)=0\}$, and the evaluation pairing descends to the quotient. Thus \cref{c:prop:arithmetic}, applied with $a=28$ to \cref{c:lem:extension}, gives an isometry $\overline U:H(L_0)\to H(L_1)$. We now choose this isometry to preserve the $k$-invariant.

\begin{lemma}\label{c:lem:k}
The isometry $\overline U$ can be chosen so that
\[
 \overline U_*(\kappa_0,0)=(\kappa_1,0).
\]
\end{lemma}
\begin{proof}
For $\ell=2,7$, the local map $f_\ell:\Z_\ell G\to P_\ell$ sends $N$ to $e$.
It therefore gives the following commutative diagram, whose rows are
\eqref{c:eq:standard-extension} and \eqref{c:eq:P-extension} tensored with $\Z_\ell$:
\[
\begin{CD}
0 @>>> \Z_\ell @>{1\mapsto N}>> \Z_\ell G @>>> (L_0)_\ell @>>> 0\\
@. @V{\id}VV @V{f_\ell}VV @V{\overline f_\ell}VV @.\\
0 @>>> \Z_\ell @>{1\mapsto e}>> P_\ell @>>> (L_1)_\ell @>>> 0.
\end{CD}
\]
The connecting maps
\[
 H^3(G;(L_i)_\ell)\xrightarrow{\cong}H^4(G;\Z_\ell)
\]
are isomorphisms, since the middle modules are projective. By \cref{c:lem:extension}, both $k$-invariants map to the same class $u$. Hence the induced map $\overline f_\ell$ sends $\kappa_0$ to $\kappa_1$, and $H(\overline f_\ell)$ sends $(\kappa_0,0)$ to $(\kappa_1,0)$.

The congruence in \cref{c:prop:arithmetic} descends under orthogonal reduction:
\[
 \overline U_\ell-H(\overline f_\ell)=28T_\ell
\]
for a $\Z_\ell G$-homomorphism $T_\ell:H((L_0)_\ell)\to H((L_1)_\ell)$. To check divisibility before passing to the quotient, both maps fix $(N,0)$ and send its orthogonal complement to $(e,0)^\perp$. If their difference on that complement is $28z$, then $\langle z,e\rangle=0$, since $\Z_\ell$ is torsion-free. The divided map therefore descends as claimed.

Positive-dimensional group cohomology is annihilated by $|G|=28$, so $\overline U_\ell$ and $H(\overline f_\ell)$ induce the same map on $H^3$. Finally, cohomology commutes with the flat extension $\Z\to\Z_\ell$: use the bar resolution, whose terms are finitely generated free modules in each degree. The finite group $H^3(G;H(L_1))$ is annihilated by $28$, and equality in its $2$- and $7$-primary parts implies equality integrally.
\end{proof}

\Needspace{12\baselineskip}
\begin{proof}[Proof of \cref{thm:noninjective}]
Mannan--Popiel's theorem gives $X\not\simeq Y$, even allowing arbitrary changes of marking. Let $M_0=M(X)$ and $M_1=M(Y)$. By \cref{c:lem:hyperbolic}, the quadratic $2$-type of $M_i$ is
\[
 \bigl(G,\hmod(L_i),(\kappa_i,0),H(L_i)\bigr).
\]
The isometry of quadratic $2$-types follows from \cref{c:prop:arithmetic,c:lem:k}. The group $Q_{28}$ has periodic cohomology of period four. Hambleton--Kreck's theorem \cite[Theorem A]{HK88} therefore gives an orientation-preserving homotopy equivalence $M_0\simeq M_1$.
\end{proof}

\bibliographystyle{amsalpha}
\bibliography{biblio}

\providecommand{\bysame}{\leavevmode\hbox to3em{\hrulefill}\thinspace}
\providecommand{\MR}{\relax\ifhmode\unskip\space\fi MR }
\providecommand{\MRhref}[2]{%
  \href{http://www.ams.org/mathscinet-getitem?mr=#1}{#2}
}
\providecommand{\href}[2]{#2}
\begin{thebibliography}{BFFM14}

\bibitem[Alp62]{Al62}
Jonathan~L. Alperin, \emph{Groups with finitely many automorphisms}, Pacific J. Math. \textbf{12} (1962), no.~1, 1--5.

\bibitem[BFFM14]{BFM14}
Eva Bayer-Fluckiger, Uriya~A. First, and Daniel~A. Moldovan, \emph{Hermitian categories, extension of scalars and systems of sesquilinear forms}, Pacific J. Math. \textbf{270} (2014), no.~1, 1--26.

\bibitem[BRK26]{BKR26}
R.~{\.I}nan{\c{c}} Baykur, Daniel Ruberman, and Robion Kirby (eds.), \emph{{K3}: A new problem list in low-dimensional topology}, Mathematical Surveys and Monographs, vol. 295, American Mathematical Society, Providence, RI, 2026.

\bibitem[Bro79]{Bro79}
Wesley~J. Browning, \emph{Finite {CW}-complexes of cohomological dimension 2 with finite abelian $\pi_1$}, Unpublished preprint, Forschungsinstitut f{\"u}r Mathematik, ETH Z{\"u}rich, 23 pp., May 1979.

\bibitem[Bro82]{Br82}
Kenneth~S. Brown, \emph{Cohomology of groups}, Graduate Texts in Mathematics, vol.~87, Springer-Verlag, New York--Berlin, 1982.

\bibitem[Car72]{Ca72}
Roger~W. Carter, \emph{Conjugacy classes in the {Weyl} group}, Compositio Math. \textbf{25} (1972), no.~1, 1--59.

\bibitem[CCP23]{CCP23}
Anthony Conway, Diarmuid Crowley, and Mark Powell, \emph{Infinite homotopy stable class for {$4$}-manifolds with boundary}, Pacific J. Math. \textbf{325} (2023), no.~2, 209--237.

\bibitem[CCPS24]{CCPS24}
Anthony Conway, Diarmuid Crowley, Mark Powell, and Joerg Sixt, \emph{Stably diffeomorphic manifolds and the realisation of modified surgery obstructions}, 2024, Version 4, 23 July 2024; to appear in M{\'e}m. Soc. Math. Fr. (N.S.).

\bibitem[CE56]{CE56}
Henri Cartan and Samuel Eilenberg, \emph{Homological algebra}, Princeton Mathematical Series, vol.~19, Princeton University Press, Princeton, New Jersey, 1956.

\bibitem[CR81]{CR81}
Charles~W. Curtis and Irving Reiner, \emph{Methods of representation theory. {V}ol. {I}}, Pure and Applied Mathematics, John Wiley \& Sons, Inc., New York, 1981, With applications to finite groups and orders. A Wiley-Interscience Publication.

\bibitem[CR87]{CR87}
\bysame, \emph{Methods of representation theory. {V}ol. {II}}, Pure and Applied Mathematics, John Wiley \& Sons, Inc., New York, 1987, With applications to finite groups and orders. A Wiley-Interscience Publication.

\bibitem[DS73a]{DySi73}
Micheal~N. Dyer and Allan~J. Sieradski, \emph{A homotopy classification of 2-complexes with finite cyclic fundamental group}, Bull. Amer. Math. Soc. \textbf{79} (1973), no.~1, 75--77.

\bibitem[DS73b]{DySi73full}
\bysame, \emph{Trees of homotopy types of two-dimensional {CW}-complexes}, Comment. Math. Helv. \textbf{48} (1973), 31--44.

\bibitem[Fox53]{Fox53}
Ralph~H. Fox, \emph{Free differential calculus. {I}. {D}erivation in the free group ring}, Ann. of Math. (2) \textbf{57} (1953), no.~3, 547--560.

\bibitem[HK88]{HK88}
Ian Hambleton and Matthias Kreck, \emph{On the classification of topological {$4$}-manifolds with finite fundamental group}, Math. Ann. \textbf{280} (1988), 85--104.

\bibitem[HK93]{HK93}
\bysame, \emph{Cancellation of hyperbolic forms and topological four-manifolds}, J. Reine Angew. Math. \textbf{443} (1993), 21--47, Preprint version: MPI/91-37.

\bibitem[HN26]{HN26}
Ian Hambleton and John Nicholson, \emph{Four-manifolds, two-complexes and the quadratic bias invariant}, Adv. Math. \textbf{486} (2026), 110728, 69 pp.; see also arXiv:2412.15089v3.

\bibitem[HN27]{HN27}
Tommy Hofmann and John Nicholson, \emph{Exotic presentations of quaternion groups and {Wall}'s {$D2$} problem}, J. Algebra \textbf{712} (2027), 253--310.

\bibitem[Jac68]{Jac68}
Heinz Jacobinski, \emph{Genera and decompositions of lattices over orders}, Acta Math. \textbf{121} (1968), 1--29.

\bibitem[Kir16]{Kir16}
Markus Kirschmer, \emph{Definite quadratic and hermitian forms with small class number}, 2016, Habilitationsschrift, RWTH Aachen University.

\bibitem[KMN22]{KMN22}
Yuya Koda, Bruno Martelli, and Hironobu Naoe, \emph{Four-manifolds with shadow-complexity one}, Ann. Fac. Sci. Toulouse Math. (6) \textbf{31} (2022), no.~4, 1111--1212.

\bibitem[Kne66]{Kn66}
Martin Kneser, \emph{Strong approximation}, Algebraic Groups and Discontinuous Subgroups, Proc. Sympos. Pure Math., vol.~9, American Mathematical Society, Providence, RI, 1966, pp.~187--196.

\bibitem[KPR22]{KPR22}
Daniel Kasprowski, Mark Powell, and Arunima Ray, \emph{Counterexamples in {$4$}-manifold topology}, EMS Surv. Math. Sci. \textbf{9} (2022), no.~1, 193--249.

\bibitem[KS84]{KS84}
Matthias Kreck and James~A. Schafer, \emph{Classification and stable classification of manifolds: some examples}, Comment. Math. Helv. \textbf{59} (1984), no.~1, 12--38.

\bibitem[Lat93]{La93}
M.~Paul Latiolais, \emph{Homotopy and homology classification of {$2$}-complexes}, Two-dimensional homotopy and combinatorial group theory (Cynthia Hog-Angeloni, Wolfgang Metzler, and Allan~J. Sieradski, eds.), London Math. Soc. Lecture Note Ser., vol. 197, Cambridge University Press, Cambridge, 1993, pp.~97--124.

\bibitem[Mar11]{Ma11}
Bruno Martelli, \emph{Four-manifolds with shadow-complexity zero}, Int. Math. Res. Not. IMRN (2011), no.~6, 1268--1351.

\bibitem[Mat89]{Mat89}
Hideyuki Matsumura, \emph{Commutative ring theory}, Cambridge Studies in Advanced Mathematics, vol.~8, Cambridge University Press, Cambridge, 1989, Translated from the Japanese by Miles Reid.

\bibitem[Met76]{Me76}
Wolfgang Metzler, \emph{{\"U}ber den {H}omotopietyp zweidimensionaler {CW}-{K}omplexe und {E}lementartransformationen bei {D}arstellungen von {G}ruppen durch {E}rzeugende und definierende {R}elationen}, J. Reine Angew. Math. \textbf{285} (1976), 7--23.

\bibitem[MP21]{MP21}
Wajid~H. Mannan and Tomasz Popiel, \emph{An exotic presentation of {$Q_{28}$}}, Algebr. Geom. Topol. \textbf{21} (2021), no.~4, 2065--2084.

\bibitem[Nic24]{Ni24}
John Nicholson, \emph{Cancellation for {$(G,n)$}-complexes and the {Swan} finiteness obstruction}, Int. Math. Res. Not. IMRN (2024), no.~16, 11832--11864.

\bibitem[Nic26]{Ni26}
\bysame, \emph{Stably free modules and the unstable classification of {$2$}-complexes}, Math. Ann. \textbf{394} (2026), Paper No. 66.

\bibitem[Phi67]{Ph67}
Anthony Phillips, \emph{Submersions of open manifolds}, Topology \textbf{6} (1967), 171--206.

\bibitem[SD79]{SiDy79}
Allan~J. Sieradski and Micheal~N. Dyer, \emph{Distinguishing arithmetic for certain stably isomorphic modules}, J. Pure Appl. Algebra \textbf{15} (1979), no.~2, 199--217.

\bibitem[Spr74]{Sp74}
Tonny~A. Springer, \emph{Regular elements of finite reflection groups}, Invent. Math. \textbf{25} (1974), 159--198.

\bibitem[Suc88]{Suc88}
Alexander~I. Suciu, \emph{The oriented homotopy type of spun {$3$}-manifolds}, Pacific J. Math. \textbf{131} (1988), no.~2, 393--399.

\bibitem[Swa60a]{Sw60}
Richard~G. Swan, \emph{Induced representations and projective modules}, Ann. of Math. (2) \textbf{71} (1960), no.~3, 552--578.

\bibitem[Swa60b]{Sw60p}
\bysame, \emph{The {$p$}-period of a finite group}, Illinois J. Math. \textbf{4} (1960), no.~3, 341--346.

\bibitem[Swa88]{Sw88}
\bysame, \emph{Torsion free cancellation over orders}, Illinois J. Math. \textbf{32} (1988), no.~3, 329--360.

\bibitem[Voi21]{Vo21}
John Voight, \emph{Quaternion algebras}, Graduate Texts in Mathematics, vol. 288, Springer, Cham, 2021.

\bibitem[Wei13]{We13}
Charles~A. Weibel, \emph{The {$K$}-book: An introduction to algebraic {$K$}-theory}, Graduate Studies in Mathematics, vol. 145, American Mathematical Society, Providence, RI, 2013.

\end{thebibliography}
\end{document}